\documentclass{amsart}
\usepackage{amssymb,latexsym}
\usepackage{amsfonts}
\usepackage{ulem}

\newtheorem{thm}{Theorem}[section]
\newtheorem{lem}[thm]{Lemma}

\newtheorem{rem}[thm]{Remark}
\newtheorem{pro}[thm]{Proposition}
\newtheorem{defi}[thm]{Definition}

\newcommand{\ba}{\begin{array}}
\newcommand{\ea}{\end{array}}

\def \qed{\cqfd}

\newcommand*{\QEDB}{\hfill\ensuremath{\square}}
\def\qed{\vbox{\hrule
\hbox{\vrule\hbox to 5pt{\vbox to 8pt{\vfil}\hfil}\vrule}\hrule}}

\newcommand{\beg}{\begin{eqnarray*}}
\newcommand{\begn}{\begin{eqnarray}}
\newcommand{\en}{\end{eqnarray*}}
\newcommand{\enn}{\end{eqnarray}}

\newcommand{\tr}{\mbox{\rm tr\,}}

\begin{document}
\title[The conical complex Monge-Amp\`ere equations on K\"ahler manifolds]{The conical complex Monge-Amp\`ere equations on K\"ahler manifolds}
\author{JiaWei Liu}
\address{Jiawei Liu\\Beijing International Center for Mathematical Research\\ Peking University \\ Beijing 100871\\ China\\} \email{jwliu@math.pku.edu.cn}
\author{Chuanjing Zhang}
\address{Chuanjing Zhang\\School of Mathematical Sciences\\
University of Science and Technology of China\\
Hefei, 230026, P.R. China\\ } \email{chjzhang@mail.ustc.edu.cn}
\thanks{}

\begin{abstract} In this paper, by providing the uniform gradient estimates for a sequence of the approximating equations, we prove the existence, uniqueness and regularity of the conical parabolic complex Monge-Amp\`ere equation with weak initial data. As an application, we prove a regularity estimates, that is, any $L^{\infty}$-solution of the conical complex Monge-Amp\`ere equation admits the $C^{2,\alpha,\beta}$-regularity.
\end{abstract}

\maketitle
\section{Introduction}
\setcounter{equation}{0}

The complex Monge-Amp\`ere equation plays an important role in geometric analysis.  It is well known that the existence of the K\"ahler-Einstein metrics and the solutions of the complex Monge-Amp\`ere equations are closely connected. Aubin \cite{TA} and Yau \cite{STY} solved the Calabi's \cite{ECA} conjecture by using the complex  Monge-Amp\`ere equation and Cao \cite{HDC} gave a parabolic proof of this conjecture by using the parabolic complex Monge-Amp\`ere equation which is equivalent to the K\"ahler-Ricci flow. By using the complex Monge-Amp\`ere equation with conical singularity,  Tian \cite{T1} and Chen-Donaldson-Sun \cite{CDS1, CDS2, CDS3} proved the Yau-Tian-Donaldson's conjecture. There are many results of the (parabolic) complex Monge-Amp\`ere equation, we refer the reader to references Caffarelli-Kohn-Nirenberg-Spruck \cite{CHNS}, J.C. Chu \cite{JCCH}, Dienw-Zhang-Zhang \cite{SDXZXWZ}, B. Guan \cite{BG2000}, Guan-Li \cite{BGQL1, BGQL2}, P.F. Guan \cite{PFG1, PFG2, PFG3}, Liu-Zhang \cite{JWLXZ}, Phong and Sturm et al. \cite{PSS, PSSW, DHPSS, PS}, G. S\'zekelyhidi \cite{GS15}, S\'zekelyhidi-Tosatti \cite{GSVT}, G. Tian \cite{T0, Tian5}, Tian-Zhu \cite{TZ}, Y. Wang \cite{YW2012}, Y.Q. Wang \cite{YQW} and X. Zhang \cite{XZ2010}, etc.

Let $(M, \omega_0)$ be a K\"ahler manifold. Assume that $D$ is an irreducible divisor, $s$ and $h$ are the definition section and smooth Hermitian metric of the line bundle associated to the divisor $D$ respectively. Denote the model conical K\"ahler metric $\omega_{\beta}=\omega_{0}+\sqrt{-1}\delta\partial\bar{\partial}|s|_{h}^{2\beta}$ and
\begin{eqnarray}
\mathcal{E}(M,\omega_{0})=\big\{\varphi\in PSH(M,\omega_{0}) | \int_{M} (\omega_{0}+\sqrt{-1}\partial\bar{\partial}\varphi)^{n}= \int_{M} \omega^{n}_{0}\big\}.
\end{eqnarray}
For any smooth function $F: \mathbb{R}\times M\rightarrow \mathbb{R}$, we consider the conical parabolic complex Monge-Amp\`ere equation
\begin{equation}\label{1}
\left \{\begin{split}
&\frac{\partial \varphi(t)}{\partial t}=\log\frac{(\omega_0+\sqrt{-1}\partial\bar{\partial}\varphi(t))^{n}}{\omega_0^{n}}+F(\varphi(t), z)+
  (1-\beta)\log|s|_h^2,\\
&\varphi(t)|_{t=0}=\varphi_0,\\
\end{split}
\right.
\end{equation}
where $\beta\in(0,1)$, $\varphi_0\in\mathcal{E}_{p}(M,\omega_{0})$ with $p>1$ and
\begin{eqnarray}
\mathcal{E}_{p}(M,\omega_{0})=\big\{\varphi\in\mathcal{E}(M,\omega_{0}) | \frac{(\omega_{0}+\sqrt{-1}\partial\bar{\partial}\varphi)^{n}}{\omega_{0}^{n}}\in L^{p}(M,\omega_{0}^{n}) \big\}.
\end{eqnarray}

\begin{defi}\label{D1}
We call $\varphi(t)$ a solution of the equation $(\ref{1})$ on $[0,T]$ $(0<T<\infty)$ if it satisfies the following conditions.
\begin{itemize}
  \item  $\varphi(t)\in C^0([0, T]\times M) \cap C^\infty((0, T]\times (M\setminus D))$  satisfies the parabolic complex Monge-Amp\`ere equation $(\ref{1})$ on $(0, T]\times (M\setminus D)$;
  \item  For any $0< \delta <T$, there exists a constant $C$ such that
   \begin{equation*}
  C^{-1}\omega_\beta \leq \omega_0+\sqrt{-1}\partial\bar{\partial}\varphi(t) \leq C \omega_\beta, \qquad on\ [\delta, T]\times (M\setminus D);
  \end{equation*}
  \item On $[\delta, T]$, there exist constants $\alpha\in(0,1)$ and $C^{\ast}$ such that $\varphi(t)$ is $C^{\alpha}$ on $M$ with respect to $\omega_{0}$ and $\| \frac{\partial\varphi(t)}{\partial t}\|_{L^{\infty}(M\setminus D)}\leqslant C^{\ast}$.
  \end{itemize}
\end{defi}

We study the conical parabolic complex Monge-Amp\`ere equation (\ref{1}) by using the following parabolic complex Monge-Amp\`ere equation
\begin{equation}\label{3}
\left \{\begin{split}
&\frac{\partial \varphi_{\varepsilon}(t)}{\partial t}=\log\frac{(\omega_0+\sqrt{-1}\partial\bar{\partial}\varphi_{\varepsilon}(t))^{n}}{\omega_{0}^{n}}+F(\varphi_{\varepsilon}(t), z)+(1-\beta)\log(\varepsilon^2+|s|_h^2),\\
&\varphi_{\varepsilon}(0)=\varphi_{0}.\\
\end{split}
\right.
\end{equation}

We first prove the uniform gradient estimates for a sequence of smooth parabolic complex Monge-Amp\`ere equations. In this process, we need construct a new auxiliary function, and we also need more details to deal with the terms which need not  be considered for a single equation, such as the terms coming from the new auxiliary function and the curvature terms, etc. Then we study the existence of the conical parabolic complex Monge-Amp\`ere equation (\ref{1}) by the approximating method which was used in \cite{JWLXZ1}. We also prove the uniqueness and regularity of the equation (\ref{1}) by using Jeffres' trick \cite{TJEF} and Tian's elliptic $C^{2,\alpha,\beta}$-estimates \cite{Tian5} respectively. In fact, we obtain the following theorem.

\begin{thm}\label{thm04} Let $(M,\omega_{0})$ be a K\"ahler manifold with complex dimension $n$, $D$ be an irreducible divisor, $s$ and $h$ be the definition section and smooth Hermitian metric of the line bundle associated to the divisor $D$ respectively. Assume that $\varphi_0\in\mathcal{E}_{p}(M,\omega_{0})$ with $p>1$. For any $\beta\in(0,1)$, there exists $T$ such that the conical parabolic complex Monge-Amp\`ere equation $(\ref{1})$ admits a unique solution $\varphi(t)$ on $[0,T]$. Furthermore, $\varphi(t)$ is $C^{2,\alpha,\beta}$ for any $\alpha\in(0,\min\{1,\frac{1}{\beta}-1\})$ when $t\in(0,T]$.
\end{thm}

\begin{rem}\label{rem05}
When $F(\varphi,z)=\mu\varphi+h$ and $c_1(M)=\mu[\omega_0]+(1-\beta)c_1(D)$ with some function $h$, the equation $(\ref{1})$ is the conical K\"ahler-Ricci flow which was studied in \cite{CW, CW1, GEDWA, JWLXZ, JWLXZ1, LMSH1, LMSH2}. When $\beta=1$, there exists no singular terms in $(\ref{1})$. Song-Tian \cite{JSGT} studied $F(z)$ case and S\'zekelyhidi-Tosatti \cite{GSVT} studied the case with particular initial data. Here, we consider $\beta\in(0,1)$ and the general initial data. 
\end{rem}

As an application of Theorem $\ref{thm04}$, we obtain the following regularity estimates for the singular equation.

\begin{thm}\label{thm06}  With the assumptions in Theorem $\ref{thm04}$, let $\varphi_{0}\in PSH(M,\omega_{0})\bigcap L^{\infty}(M)$ be a solution of the conical complex Monge-Amp\`ere equation \begin{equation}\label{2016080501}
(\omega_{0}+\sqrt{-1}\partial\bar{\partial}\varphi_{0})^{n}=e^{-F(\varphi_{0}, z)}\frac{\omega_{0}^{n}}{|s|_{h}^{2(1-\beta)}}.
\end{equation}
Then  $\varphi_{0}\in C^{2,\alpha,\beta}(M)$.
\end{thm}

\begin{rem}\label{rem06}  In \cite{GSVT}, S\'zekelyhidi-Tosatti considered this problem in the smooth case $\beta=1$. In their case, they can prove that $\varphi_{0}$ is smooth. In \cite{Tian5}, Tian obtained the $C^{2,\alpha,\beta}$-regularity for $\varphi_{0}$ under the assumption that $\omega_{0}+\sqrt{-1}\partial\bar{\partial}\varphi_{0}$ is equivalent to the model conical K\"ahler metric $\omega_{\beta}$. Here we only assume that  $\varphi_{0}\in PSH(M,\omega_{0})\bigcap L^{\infty}(M)$. When $F(\varphi,z)=\mu\varphi+h$ with smooth function $h$, this result is proved by Guenancia-P$\breve{a}$un in \cite{GP1}. In our case, we need the uniform gradient estimates which do not be needed in Guenancia-P$\breve{a}$un's case for the sequence of  approximating equations. \end{rem}

The paper is organized as follows. In section $2$, we prove the uniform gradient estimates for the sequence of the approximating equations, and then we prove the existence and uniqueness of the parabolic complex Monge-Amp\`ere equation $(\ref{3})$. In section $3$,  we prove the existence of the solution to the conical parabolic complex Monge-Amp\`ere equation $(\ref{1})$ by limiting the equations $(\ref{3})$,  and prove that $\varphi(t)$ converges to $\varphi_{0}$ in $L^{\infty}$-norm as $t\rightarrow0^{+}$. We also prove the uniqueness and regularity of the conical equation $(\ref{1})$. At last, we prove Theorem $\ref{thm06}$.

{\bf  Acknowledgement:} The  authors would like to thank their advisors Professor Jiayu Li and Professor Xi Zhang for providing many suggestions and encouragements. The first author also would like to thank Professor Xiaohua Zhu  for his constant help.

\section{The existence of the parabolic complex Monge-Amp\`ere equation with weak initial data}
\setcounter{equation}{0}

In this section, we prove the existence and uniqueness of the solution to the parabolic complex Monge-Amp\`ere equation $(\ref{3})$ by following the arguments in \cite{JWLXZ1, JSGT}.  The differences are that we need prove the uniform gradient estimates for the sequence of the approximating equations in this case. In \cite{GSVT}, Sz\'ekelyhidi and Tosatti proved the gradient estimate for a single parabolic Monge-Amp\`ere equation. But when we prove the uniform gradient estimates for a sequence of parabolic Monge-Amp\`ere equations, we need more details to deal with the terms which need not be considered for a single equation.

We denote $f=\frac{(\omega_{0}+\sqrt{-1}\partial\bar{\partial}\varphi_{0})^{n}}{\omega_{0}^{n}}\in L^{p}(M,\omega_{0}^{n})$ for some $p>1$. By considering the complex Monge-Amp\`ere equation
\begin{eqnarray}
(\omega_{0}+\sqrt{-1}\partial\bar{\partial}\varphi_{0,j})^{n}=f_{j}\omega_{0}^{n}
\end{eqnarray}
and using the stability theorem in \cite{K00} (see also \cite{SDZZ} or \cite{VGAZ1}), we have
\begin{eqnarray}\label{00000}
\lim\limits_{j\rightarrow\infty}\|\varphi_{0,j}-\varphi_{0}\|_{L^{\infty}(M)}=0,
\end{eqnarray}
where $\varphi_{0,j}\in PSH(M,\omega_0)\cap C^\infty(M)$ satisfy
$\sup\limits_{M}(\varphi_0-\varphi_{0,j})=\sup\limits_{M}(\varphi_{0,j}-\varphi_{0})$, $f_{j}\in C^{\infty}(M)$ satisfy $\int_{M}f_{j}\omega_{0}^{n}=\int_{M}\omega_{0}^{n}$ and $
\lim\limits_{j\rightarrow\infty}\|f_{j}-f\|_{L^{p}(M)}=0$.

We prove the existence of the solution to the equation $(\ref{3})$ by using the smooth parabolic complex Monge-Amp\`ere equation
\begin{equation}\label{CMAE1}
\left \{\begin{split}
&\frac{\partial \varphi_{\varepsilon,j}(t)}{\partial t}=\log\frac{(\omega_0+\sqrt{-1}\partial\bar{\partial}\varphi_{\varepsilon,j}(t))^{n}}{\omega_{0}^{n}}+F(\varphi_{\varepsilon, j}(t), z)+(1-\beta)\log(\varepsilon^2+|s|_h^2),\\
&\varphi_{\varepsilon,j}(0)=\varphi_{0,j},\\
\end{split}
\right.
\end{equation}
which can be written as
\begin{equation}\label{CMAE2}
\left \{\begin{split}
\frac{\partial \phi_{\varepsilon, j}(t)}{\partial t}=&\log\frac{(\omega_\varepsilon+\sqrt{-1}\partial\bar{\partial}\phi_{\varepsilon,j}(t))^{n}}{\omega_{\varepsilon}^{n}}+F(\phi_{\varepsilon, j}(t)+\delta\chi(\varepsilon^2+|s|_h^2), z)+f_{\varepsilon},\\
\phi_{\varepsilon,j}(0)=&\varphi_{0,j}-\delta\chi(\varepsilon^2+|s|_h^2):=\phi_{\varepsilon,0,j},\\
\end{split}
\right.
\end{equation}
where $\phi_{\varepsilon,j}(t)=\varphi_{\varepsilon,j}(t)-\delta\chi(\varepsilon^2+|s|_h^2)$, $\chi(\varepsilon^{2}+|s|_h^2)=\frac{1}{\beta}\int_{0}^{|s|_h^2}\frac{(\varepsilon^{2}+r)^{\beta}-
\varepsilon^{2\beta}}{r}dr$, $\omega_{\varepsilon}=\omega_{0}+\sqrt{-1}\delta\partial\overline{\partial}\chi(\varepsilon^{2}+|s|_{h}^{2})$, $f_{\varepsilon}=\log(\frac{\omega_{\varepsilon}^{n}}{\omega_{0}^{n}}\cdot(\varepsilon^{2}+|s|_{h}^{2})^{1-\beta})$. We know that $\chi(\varepsilon^2+|s|_h^2)$ and $f_{\varepsilon}$ are uniformly bounded (see $(15)$ and $(25)$ in \cite{CGP}).

Let $[0,T_{\varepsilon,j})$ be the maximal existence interval of the equation $(\ref{CMAE2})$, where $T_{\varepsilon,j}$ depends on the $C^{2,\alpha}$-norm of $\phi_{\varepsilon,0,j}$. Our aim is to obtain the uniform high order estimates, which only depend on the initial condition in a weaker way. For all $s$, we define
\begin{equation}
\overline{F}(s)=\max_{\substack{\varepsilon\in [0, 1]\\ z \in M}}F(s+\delta\chi(\varepsilon^2+|s|_h^2), z)+C,
\end{equation}
which is a continuous function independent of $\varepsilon$ and $j$, where $C$ satisfies $|f_{\varepsilon}|\leq C$. Let $G(t)$ be the solution of equation
\begin{equation}
\left \{\begin{split}
&G'(t)=\overline{F}(G(t))\\
&G(0)=A\\
\end{split}
\right.
\end{equation}
on $[0,\overline{T}]$, where $ \|\phi_{\varepsilon, 0, j}\|_{L^\infty(M)}\leqslant A$ for any $\varepsilon$ and $j$. Denote $\overline{T}_{\varepsilon,j}=\min\limits\{T_{\varepsilon,j},\overline{T}\}$.

\begin{lem}\label{2016080701}
For any $0<T<\overline{T}_{\varepsilon,j}$, there exists a constant $C$ depending only on $\|\varphi_0\|_{L^\infty(M)}$, $\beta$, $n$, $\omega_{0}$, $\overline{T}$ and $F$ such that for any $t\in[0,T]$ and $\varepsilon>0$,
\begin{equation}
\|\phi_{\varepsilon, j}(t)\|_{L^\infty(M)}\leq C.
\end{equation}
\end{lem}
{\bf Proof. }  For any $\varepsilon>0$, let  $z_{t}$ be the maximum point of $\phi_{\varepsilon,j}$ at time $t$.
\begin{equation}
\begin{split}
\frac{\partial \phi_{\varepsilon, j}}{\partial t}(t, z_t)=&\log\frac{(\omega_0+\sqrt{-1}\partial\bar{\partial}\phi_{\varepsilon,j}(t, z_t))^{n}}{\omega_{\varepsilon}^{n}}\\
&+F(\phi_{\varepsilon, j}(t, z_t)+\delta\chi(\varepsilon^2+|s|_h^2), z_t)+f_{\varepsilon}(z_t)\\
\leq& \overline{F}(\phi_{\varepsilon, j}(t, z_t)).\\
\end{split}
\end{equation}
If the derivative does not exist at some point, we consider the limsup of the forward difference quotients at these points ( Hamilton \cite{RSHA}). By the maximum principle (see Morgan-Tian \cite{JWMGT}), we have $\max\limits_M \phi_{\varepsilon, j}(t)\leq G(t)$, where the choice of function $G$ is independent of $\varepsilon$ and $j$. Hence there exists a constant $C$ depending only on $\|\varphi_0\|_{L^\infty(M)}$, $\beta$, $n$, $\omega_{0}$, $\overline{T}$ and $F$ such that $\max\limits_{M\times[0,T]} \phi_{\varepsilon, j}(t)\leq C$ for any $0<T<\overline{T}_{\varepsilon,j}$, $\varepsilon>0$ and $j$.
By the similar arguments we get the lower bound. \QEDB

\begin{pro}
 For any $0<T< \overline{T}_{\varepsilon,j}$, there exists a constant $C$ depending only on $\|\varphi_0\|_{L^\infty(M)}$, $n$, $\beta$, $F$, $\overline{T}$ and $\omega_{0}$ such that for any $t\in(0,T]$ and $\varepsilon>0$,
\begin{equation}\label{2011111}
\frac{t^n}{C}\leq\frac{(\omega_{\varepsilon}+\sqrt{-1}\partial\bar{\partial}\phi_{\varepsilon,j}(t))^{n}}{\omega_{\varepsilon}^n}\leq e^{\frac{C}{t}}.
\end{equation}
\end{pro}
{\bf Proof:} Let $\Delta_{\varepsilon,j}$ be the Laplacian operator associated to the K\"ahler form $\omega_{\varepsilon,j}(t)=\omega_{\varepsilon}+\sqrt{-1}\partial\bar{\partial}\phi_{\varepsilon,j}(t)$. Straightforward calculations show that
\begin{equation}
(\frac{\partial}{\partial t}-\Delta_{\varepsilon, j})\dot{\phi}_{\varepsilon,j}=F'(\phi_{\varepsilon, j}(t)+k\chi(\varepsilon^2+|s|_h^2), z)\dot{\phi}_{\varepsilon,j}.
\end{equation}
Let $H_{\varepsilon,j}^+(t)=t\dot{\phi}_{\varepsilon,j}(t)-A\phi_{\varepsilon,j}(t)$. Then
\begin{equation}
\begin{split}
(\frac{\partial}{\partial t}-\Delta_{\varepsilon,j})H_{\varepsilon,j}^+(t)=(1+F't-A)\dot{\phi}_{\varepsilon,j}+An-A\tr_{\omega_{\varepsilon,j}}\omega_\varepsilon.
\end{split}
\end{equation}
Let $(t_0,x_0)$ be the maximum point of $H_{\varepsilon,j}^+(t)$ on $[0,T]\times M$. Without loss of generality, we can assume that $t_{0}>0$ and $\dot{\phi}_{\varepsilon,j}(t_0, x_0)> 0$. By the maximum principle,
\begin{equation}
\dot{\phi}_{\varepsilon,j}(t)\leq \frac{C}{t}\ \ on\ \ (0,T]\times M.
\end{equation}
Let $H_{\varepsilon,j}^-(t)=\dot{\phi}_{\varepsilon,j}(t)+B\phi_{\varepsilon,j}(t)-n\log t$. Then $H^-_{\varepsilon,j}(t)$ tends to $+\infty$ as $t\rightarrow0^+$ and
\begin{equation}
\begin{split}
(\frac{\partial}{\partial t}-\Delta_{\varepsilon,j})H_{\varepsilon,j}^-(t)=(F'+B)\dot{\phi}_{\varepsilon,j}+B\tr_{\omega_{\varepsilon,j}}\omega_\varepsilon-Bn-\frac{n}{t}.\end{split}
\end{equation}
Assume that $(t_0,x_0)$ is the minimum point of $H_{\varepsilon,j}^-(t)$ on $[0,T]\times M$. We conclude that $t_0>0$ and there exists constants $C_1$, $C_2$ and $C_3$ such that, at $(t_{0},x_{0})$,
\begin{eqnarray}\label{2002}\nonumber(\frac{\partial}{\partial t}-\Delta_{\varepsilon,j})H_{\varepsilon,j}^-(t)|_{(t_0,x_0)}&\geq& \big(C_1\big(\frac{\omega_{\varepsilon}^n}{\omega^{n}_{\varepsilon,j}(t)}\big)^{\frac{1}{n}}+C_2\log\frac{\omega^n_{\varepsilon,j}(t)}{\omega_{\varepsilon}^n}
-\frac{C_3}{t}\big)|_{(t_0,x_0)}\\
&\geq&\big(\frac{C_1}{2}\big(\frac{\omega_{\varepsilon}^n}{\omega^{n}_{\varepsilon,j}(t)}\big)^{\frac{1}{n}}
-\frac{C_3}{t}\big)|_{(t_0,x_0)},
\end{eqnarray}
where the constant $C_1$ depends only on $n$, $C_2$ depends only on $F$ and $C_3$ depends only on $n$, $\omega_{0}$, $\|\varphi_0\|_{L^\infty(M)}$, $\beta$, $F$ and $\overline{T}$. In inequality $(\ref{2002})$, without loss of generality, we assume that $\frac{\omega_{\varepsilon}^n}{\omega^{n}_{\varepsilon,j}(t)}>1$ and $\frac{C_1}{2}(\frac{\omega_{\varepsilon}^n}{\omega^{n}_{\varepsilon,j}(t)})^{\frac{1}{n}}+C_2\log\frac{\omega^n_{\varepsilon,j}(t)}
{\omega_{\varepsilon}^n}\geq0$ at $(t_0,x_0)$. By the maximum principle, we have
\begin{eqnarray}\omega^{n}_{\varepsilon,j}(t_0,x_0)\geq C_4t^n \omega_{\varepsilon}^n(x_0),
\end{eqnarray}
where $C_4$ is independent of $\varepsilon$ and $j$. Then it easily follows that $H_{\varepsilon,j}^-(t)$ is bounded from below by a constant $C$ depending only on $\|\varphi_0\|_{L^\infty(M)}$, $n$, $\beta$, $\omega_{0}$, $F$ and $\overline{T}$.\QEDB

\medskip

Next, we prove the uniform gradient estimates for equation $(\ref{CMAE2})$.
\begin{pro}\label{2016080703}
For any $0<T<\overline{T}_{\varepsilon,j}$, there exists a constant $C$ depending only on $\|\varphi_0\|_{L^\infty(M)}$, $n$, $\beta$, $F$, $\omega_{0}$ and $\overline{T}$ such that for any $t\in(0,T]$, $\varepsilon>0$ and $j$,
\begin{equation}
|\nabla \phi_{\varepsilon,j}(t)|_{\omega_{\varepsilon}}\leq e^{e^{\frac{C}{t}}}.
\end{equation}
\end{pro}

{\bf Proof.} We modify Blocki's estimates \cite{ZB1} for the complex Monge-Amp\`ere equation (compare Hanani \cite{AHAN} and Sz\'ekelyhidi-Tosatti  \cite{GSVT}). Define
\begin{equation}
H_{\varepsilon, j}(t)= e^{-\frac{\alpha}{t}}\log|\nabla \phi_{\varepsilon, j}|_{\omega_{\varepsilon}}^2-\gamma(\phi_{\varepsilon, j})+At\log t+B \Psi_{\varepsilon\rho}e^{-\frac{\alpha}{t}},
\end{equation}
where function $\gamma$, constants $\alpha$, $A$ and $B$ will be chosen later and $\Psi_{\varepsilon\rho}=\frac{1}{\rho}\int_{0}^{|s|_h^2}\frac{(\varepsilon^{2}+r)^{\rho}-\varepsilon^{2\rho}}{r}dr$. Let $(t_{0}, x_{0})$ be the maximum point of $H_{\varepsilon,j}(t)$ on $[0,T]\times M$, we need only consider $t_0>0$. We choose a local coordinate system $w = (w_{1}, . . . , w_{n})$, to make $(g_{\varepsilon i\bar{j}})$ be identity and $(g_{\phi_{\varepsilon,j}i\bar{k}})$ be a diagonal matrix. For simplicity, we write $\Phi=|\nabla \phi_{\varepsilon,j}(t)|^{2}_{\omega_{\varepsilon}}$ and $\phi_{\varepsilon,j}(t)=\phi$. Straightforward calculations show that
\begin{eqnarray}\label{2016080702}\nonumber
(\frac{\partial}{\partial t}-\Delta_{\varepsilon,j})H_{\varepsilon,j}(t)&=&\frac{\alpha}{t^2}e^{-\frac{\alpha}{t}}\log|\nabla \phi|_{\omega_{\varepsilon}}^2+A\log t+A-B\Delta_{\varepsilon,j}\Psi_{\varepsilon\rho}e^{-\frac{\alpha}{t}}-\gamma'\dot{\phi}\\\nonumber
&&-\sum\limits_{j,k,p}e^{-\frac{\alpha}{t}}\frac{R_{\varepsilon j\bar{k}p\bar{p}}\phi_j\phi_{\bar{k}}}{\Phi (1+\phi_{p\bar{p}})}+e^{-\frac{\alpha}{t}}\sum\limits_{j}\frac{2Re F' \phi_{j}\phi_{\bar{j}}}{\Phi}+\gamma''|\nabla\phi|_{\omega_{\varepsilon, j}}^2+\gamma'n\\
&&+e^{-\frac{\alpha}{t}}\delta\sum\limits_{j}\frac{2Re F' \chi_{j}\phi_{\bar{j}}}{\Phi}+e^{-\frac{\alpha}{t}}\sum\limits_{j}\frac{2Re f_{\varepsilon j}\phi_{\bar{j}}}{\Phi}+e^{-\frac{\alpha}{t}}\sum\limits_{j}\frac{2Re F_j\phi_{\bar{j}}}{\Phi}\\\nonumber
&&-e^{-\frac{\alpha}{t}}\sum\limits_{j,p}\frac{|\phi_{jp}|^2}{\Phi (1+\phi_{p\bar{p}})}-e^{-\frac{\alpha}{t}}\sum\limits_{p}\frac{|\phi_{p\bar{p}}|^2}{\Phi (1+\phi_{p\bar{p}})}+\frac{\alpha}{t^2}e^{-\frac{\alpha}{t}}B\Psi_{\varepsilon\rho}-\gamma'\tr_{\omega_{\varepsilon, j}}\omega_{\varepsilon}\\\nonumber
&&+\sum\limits_{p}\frac{e^{\frac{\alpha}{t}}(\gamma')^2|\phi_{p}|^2+e^{-\frac{\alpha}{t}}B^2|\Psi_{\varepsilon\rho p}|^2-2B\gamma'Re\phi_p\Psi_{\varepsilon\rho\bar{p}}}{1+\phi_{p\bar{p}}},
\end{eqnarray}
where $F'$ is the derivative of $F(\phi, z)$ in the $\phi$ variable and $F_{k}$ is the derivative  of $F(\phi, z)$ in the $z$ variable. At $(x_{0}, t_{0})$, we have
\begin{eqnarray}\sum_{k} \phi_{kp}\phi_{\bar{k}}=e^{\frac{\alpha}{t}}\Phi\phi_{p}\gamma^{\prime}-B\Phi\Psi_{\varepsilon\rho p}-\phi_{p}\phi_{p\bar{p}}.
\end{eqnarray}
It follows that
\begin{eqnarray}\label{2016040302} \nonumber e^{-\frac{\alpha}{t}}\sum_{k,p}\frac{|\phi_{kp}|^{2}}{\Phi(1+\phi_{p\bar{p}})}&\geq& e^{-\frac{\alpha}{t}}\sum_{p}\frac{|e^{\frac{\alpha}{t}}\Phi\phi_{p}\gamma^{\prime}-B\Phi\Psi_{\varepsilon\rho p}-\phi_{p}\phi_{p\bar{p}}|^{2}}{\Phi^{2}(1+\phi_{p\bar{p}})}\\
&\geq&e^{\frac{\alpha}{t}}\sum_{p}\frac{\gamma^{\prime2}|\phi_{p}|^{2}}{1+\phi_{p\bar{p}}}+e^{-\frac{\alpha}{t}}\sum_{p}\frac{B^{2}|\Psi_{\varepsilon\rho p}|^{2}}{1+\phi_{p\bar{p}}}-\sum_{p}\frac{2B\gamma^{\prime}Re\phi_{p}\Psi_{\varepsilon\rho \bar{p}}}{1+\phi_{p\bar{p}}}\\ \nonumber
&\ &-\sum\limits_{p}\frac{2\gamma'|\phi_p|^2\phi_{p\bar{p}}}{\Phi(1+\phi_{p\bar{p}})}-\sum_{p}\frac{2e^{-\frac{\alpha}{t}}B Re\phi_p\phi_{p\bar{p}}\Psi_{\varepsilon\rho\bar{p}}}{\Phi(1+\phi_{p\bar{p}})},
\end{eqnarray}
where we assume that $\gamma'>0$, $e^{-\frac{\alpha}{2t_0}}\log\Phi(x_{0},t_{0})\geq C$ for some constant $C$ which will be determined later, and hence $\Phi(x_{0},t_{0})\geq1$. By computing, we have
\begin{eqnarray}\label{3.22.23}|\nabla\chi|^{2}_{\omega_{\varepsilon}}\leq C\big(\frac{(\varepsilon^{2}+|s|^{2}_{h})^{\beta}-
\varepsilon^{2\beta}}{\beta|s|_{h}^{2}}\big)^{2}|s|^{2}_{h}(\varepsilon^{2}+|s|^{2}_{h})^{1-\beta}|D's|^{2}_{h}+C,\end{eqnarray}
where constant $C$ is uniform. Since the function $$F(x,y)=\frac{((x^{2}+y^{2})^{\beta}-x^{2\beta})^{2}}{y^{2}}(x^{2}+y^{2})^{1-\beta}$$ is bounded on $[0,1]\times[0,1]$ and $|\nabla\chi|^{2}_{\omega_{\varepsilon}}$ is uniformly bounded.
\begin{equation}
\begin{split}
&\sum\limits_{p}\Big(\frac{e^{\frac{\alpha}{t}}(\gamma')^2|\phi_{p}|^2+e^{-\frac{\alpha}{t}}B^2|\Psi_{\varepsilon\rho p}|^2-2B\gamma'Re\phi_p\Psi_{\varepsilon\rho\bar{p}}}{u_{p\bar{p}}}-e^{-\frac{\alpha}{t}}\frac{|\phi_{p\bar{p}}|^2}{\Phi (1+\phi_{p\bar{p}})}\Big)\\
\leq& \sum\limits_{p}\frac{2 B e^{-\frac{\alpha}{t}} |\phi_{p\bar{p}}||\Psi_{\varepsilon\rho\bar{p}}|}{\sqrt{\Phi}(1+\phi_{p\bar{p}})}+\sum\limits_{p}\frac{2\gamma'|\phi_p|^2\phi_{p\bar{p}}}{\Phi(1+\phi_{p\bar{p}})}\\
\leq &2\gamma'n+2B^2e^{-\frac{\alpha}{t}}\frac{|\Psi_{\varepsilon\rho\bar{p}}|^2}{1+\phi_{p\bar{p}}}+\frac{e^{-\frac{\alpha}{t}}|\phi_{p\bar{p}}|^2}{2\Phi(1+\phi_{p\bar{p}})},
\end{split}
\end{equation}
\begin{eqnarray}\label{2016040303}\nonumber
2\sum_{k}\frac{Ref_{\varepsilon k}\phi_{\bar{k}}}{\Phi}&\leq&2\sum_{k}\frac{\sum_{l}(1+\phi_{l\bar{l}})}{\sqrt{\Phi}\sqrt{\sum_{l}(1+\phi_{l\bar{l}})}}\cdot\frac{|f_{\varepsilon k}|}{\sqrt{\sum_{l}(1+\phi_{l\bar{l}})}}\\
&\leq&\frac{C}{K\Phi}\sum_{l}\frac{1+|\phi_{l\bar{l}}|^{2}}{1+\phi_{l\bar{l}}}+K\sum_{k}\frac{|f_{\varepsilon k}|^{2}}{1+\phi_{k\bar{k}}}.
\end{eqnarray}
By choosing $K$ large enough and putting $(\ref{2016040302})$-$(\ref{2016040303})$ into $(\ref{2016080702})$, we have
\begin{eqnarray}\label{2016040304}\nonumber(\frac{\partial}{\partial t}-\Delta_{\varepsilon,j})H_{\varepsilon,j}(t)&\leq&\frac{\alpha}{t^{2}}e^{-\frac{\alpha}{t}}\log\Phi+(A-C\gamma^{\prime})\log t+C\gamma^{\prime}-Be^{-\frac{\alpha}{t}}\Delta_{\varepsilon,j}\Psi_{\varepsilon\rho}\\
&\ &-e^{-\frac{\alpha}{t}}\sum_{k,l,p}\frac{R_{l\bar{k}p\bar{p}}\phi_{l}\phi_{\bar{k}}}{\Phi(1+\phi_{p\bar{p}})}+Ke^{-\frac{\alpha}{t}}\sum_{k}\frac{|f_{\varepsilon k}|^{2}}{1+\phi_{k\bar{k}}}+2B^2e^{-\frac{\alpha}{t}}\frac{|\Psi_{\varepsilon\rho\bar{p}}|^2}{1+\phi_{p\bar{p}}}\\ \nonumber
&\ &+(C-\gamma^{\prime})\sum_{k}\frac{1}{1+\phi_{k\bar{k}}}+\gamma^{\prime\prime}\sum_{k}\frac{|\phi_{k}|^{2}}{1+\phi_{k\bar{k}}}+C.
\end{eqnarray}
We claim that
\begin{eqnarray}\label{2016040305'}\nonumber&\ &-B\Delta_{\varepsilon,j}\Psi_{\varepsilon\rho}-\sum_{k,l,p}\frac{R_{l\bar{k}p\bar{p}}\phi_{l}\phi_{\bar{k}}}{\Phi(1+\phi_{p\bar{p}})}+K\sum_{k}\frac{|f_{\varepsilon k}|^{2}}{1+\phi_{k\bar{k}}}+2B^2\frac{|\Psi_{\varepsilon\rho\bar{p}}|^2}{1+\phi_{p\bar{p}}}\\\nonumber
&\leq&B\sum_{k}\frac{1}{1+\phi_{k\bar{k}}}-B\sum_{k}\frac{1}{(\varepsilon^2+|z^{n}|^2)^{1-\rho}}|\frac{\partial z^n}{\partial w^k}|^2\frac{1}{1+\phi_{k\bar{k}}}\\
&\ &+2B^2\sum_k\frac{((\varepsilon^2+|s|_h^2)^{\rho}-\varepsilon^{2\rho})^2}{|s|_h^2}\big|\frac{\partial z^n}{\partial\omega^k}\big|^{2}\frac{1}{1+\phi_{k\bar{k}}}\\\nonumber
&\ &+C\sum_{k}\frac{1}{(\varepsilon^2+|z^{n}|^2)^{\tilde{\beta}}}|\frac{\partial z^n}{\partial w^k}|^2\frac{1}{1+\phi_{k\bar{k}}},
\end{eqnarray}
where $\tilde{\beta}\in\{\beta, 1-\beta, 2\beta-1, 1-2\beta\}$. Set $\theta\in[0, \frac{\pi}{2})$, $\varepsilon=r\cos\theta$, $|s|_h=r\sin\theta$ and $r\in(0, 1)$. Then
\begin{eqnarray}\frac{((\varepsilon^2+|s|_h^2)^{\rho}-\varepsilon^{2\rho})^2}{|s|_h^2}=\frac{r^{4\rho}(1-\cos^{2\rho}\theta)^2}{r^2\sin^2\theta}.
\end{eqnarray}
It is obvious that
\begin{equation}
\frac{r^{4\rho}(1-\cos^{2\rho}\theta)^2}{r^2\sin^2\theta}\leq C_{1}\frac{r^{2\rho}}{r^{2(1-\rho)}}.
\end{equation}
We first confirm $B\geq2C$ in $(\ref{2016040305'})$. Fix a neighbourhood $U$ of divisor $D$ such that
\begin{equation}
2BC_{1}r^{2\rho}< \frac{1}{2}\ \ in\ \ U.
\end{equation}
We need only consider the case $x_{0}\in U$. By $1-\rho>\tilde{\beta}$, we have
\begin{eqnarray}\label{2016040305}\nonumber(\frac{\partial}{\partial t}-\Delta_{\varepsilon,j})H_{\varepsilon,j}(t)&\leq&e^{-\frac{\alpha}{2t}}\log\Phi+(A-C\gamma^{\prime})\log t+C\gamma^{\prime}+C\\ \nonumber
&\ &+(C-\gamma^{\prime})\sum_{k}\frac{1}{1+\phi_{k\bar{k}}}+\gamma^{\prime\prime}\sum_{k}\frac{|\phi_{k}|^{2}}{1+\phi_{k\bar{k}}}.
\end{eqnarray}
Let $\gamma(s)=D's-\frac{1}{D'}s^{2}$. If we choose $D'$ large enough and $A$ satisfying $A-C\gamma^{\prime}\geq0$,
\begin{eqnarray}\label{2016040306}\sum_{k}\frac{1}{1+\phi_{k\bar{k}}}+\sum_{k}\frac{|\phi_{k}|^{2}}{1+\phi_{k\bar{k}}}&\leq&2e^{-\frac{\alpha}{2t}}\log\Phi,
\end{eqnarray}
so in particular $(1+\phi_{k\bar{k}})^{-1}\leq 2e^{-\frac{\alpha}{2t}}\log\Phi$ for any $k$. By $(\ref{2011111})$, we know that
\begin{eqnarray}\label{2016040307}1+\phi_{k\bar{k}}\leq Ce^{\frac{C-\alpha(n-1)}{2t}}(\log\Phi)^{n-1},\ \ \ k=1,\cdots,n.
\end{eqnarray}
If we choose $\alpha$ large enough, we have
\begin{eqnarray}\label{2016040308}1+\phi_{k\bar{k}}\leq C(\log\Phi)^{n-1},\ \ \ k=1,\cdots,n.
\end{eqnarray}
By $(\ref{2016040306})$, we conclude that at $(x_{0},t_{0})$,
\begin{eqnarray}\label{2016040309}\Phi\leq C(\log\Phi)^{n},\ \ \ k=1,\cdots,n.
\end{eqnarray}
This shows that $\Phi(t_{0},x_{0})\leq C$ and in turn $H_{\varepsilon,j}(t_{0},x_{0})$ are uniform bounded by some uniform constant. Hence there exists a constant $C$ depending only on $\|\varphi_0\|_{L^\infty(M)}$, $n$, $\beta$, $\omega_{0}$, $F$ and $\overline{T}$ such that
\begin{eqnarray}|\nabla \phi_{\varepsilon,j}(t)|^{2}_{\omega_{\varepsilon}}\leq e^{e^{\frac{C}{t}}}.
\end{eqnarray}

\begin{lem}\label{205} For any $0<T<\overline{T}_{\varepsilon,j}$, there exists a constant $C$ depending only on $\|\varphi_0\|_{L^\infty(M)}$, $n$, $\beta$, $\omega_{0}$, $F$ and $\overline{T}$ such that for any $t\in(0,T]$ and $\varepsilon>0$,
\begin{eqnarray}e^{-e^{e^{\frac{C}{t}}}}\omega_\varepsilon\leq\omega_{\varepsilon,j}(t)\leq e^{e^{e^{\frac{C}{t}}}}\omega_\varepsilon.
\end{eqnarray}
\end{lem}

{\bf Proof:}\ \ Let $\omega_{\varepsilon,j}=\omega_\varepsilon+\sqrt{-1}\partial\bar{\partial}\phi_{\varepsilon,j}(t)$. By computing, we have
\begin{equation}
\begin{split}
&(\frac{d}{dt}-\Delta_{\varepsilon, j})\log\tr_{\omega_\varepsilon}\omega_{\varepsilon, j}\\
=&\frac{1}{\tr_{\omega_\varepsilon}\omega_{\varepsilon, j}}\big(\Delta_{\omega_{\varepsilon}}(\dot{\varphi}-\log\frac{\omega_{\varepsilon, j}^n}{\omega_{\varepsilon}^n})+R_{\omega_{\varepsilon}}\big)-\frac{1}{\tr_{\omega_\varepsilon}\omega_{\varepsilon, j}}(g_{\varepsilon, j}^{p\bar{q}}g_{\varepsilon,j\ m\bar{k}}R_{\omega_{\varepsilon\ p\bar{q}}}^{\bar{m}k})\\
&+\frac{g_{\varepsilon, j}^{\delta\bar{k}}\partial_{\delta}\tr_{\omega_{\varepsilon}}\omega_{\varepsilon, j}\partial_{\bar{k}}\tr_{\omega_{\varepsilon}}\omega_{\varepsilon, j}}{(\tr_{\omega_{\varepsilon}}\omega_{\varepsilon, j})^2}-\frac{g_\varepsilon^{\gamma \bar{s}}\varphi_{\gamma\ p}^{\ \ t}\varphi_{\bar{s}t}^{\ \ p}}{\tr_{\omega_{\varepsilon}}\omega_{\varepsilon, j}}.
\end{split}
\end{equation}
At the same time, by using Proposition $\ref{2016080703}$, we have
\begin{equation}
\begin{split}
\Delta_{\omega_\varepsilon}F(\phi_{\varepsilon, j}+k\chi, z)=&F'(\Delta_{\omega_\varepsilon}\phi_{\varepsilon, j}+k\Delta_{\omega_{\varepsilon}}\chi)+\Delta_{\omega_\varepsilon}F\\
&+CF''(|\nabla\phi_{\varepsilon, j}|^2_{\omega_\varepsilon}+k|\nabla\chi|_{\omega_{\varepsilon}}^2)+2|F'_k|_{\omega_\varepsilon}(|\nabla\phi_{\varepsilon, j}|_{\omega_\varepsilon}
+k|\nabla\chi|_{\omega_\varepsilon})\\
\leq &Ce^{e^{\frac{C}{t}}}+C\tr_{\omega_\varepsilon}\omega_{\varphi_\varepsilon}.
\end{split}
\end{equation}
Then by the similar arguments as that in Lemma $2.3$ \cite{JWLXZ1}, we have
\begin{equation}
\begin{split}
&(\frac{\partial}{\partial t}-\Delta_{\varepsilon, j})(e^{-e^{\frac{\alpha}{t}}}\log\tr_{\omega_{\varepsilon}}
\omega_{\varepsilon, j}+e^{-e^{\frac{\alpha}{t}}}\Psi_{\varepsilon\rho}-A\phi_{\varepsilon, j})\\
\leq & -\frac{A}{4}\tr_{\omega_{\varepsilon, j}}\omega_{\varepsilon}-C\log t+C.
\end{split}
\end{equation}
By the maximum principle and the inequality
\begin{eqnarray}\label{206}
tr_{\omega_{\varepsilon}}\omega_{\varepsilon,j}(t)\leq\frac{1}{(n-1)!}(tr_{\omega_{\varepsilon,j}(t)}\omega_{\varepsilon})^{n-1}\frac{\omega^n_{\varepsilon,j}(t)}{\omega_{\varepsilon}^n},
\end{eqnarray}
we conclude
\begin{equation}\label{207}
\tr_{\omega_{\varepsilon}}\omega_{\varepsilon, j}\leq e^{h_{C}(t)},
\end{equation}
where $h_\alpha(t)=e^{e^{\frac{\alpha}{t}}}$. Using the inequality $(\ref{206})$ again, we conclude that
\begin{equation}\label{208}
\tr_{\omega_{\varepsilon, j}}\omega_{\varepsilon}\leq e^{h_{C}(t)},
\end{equation}
where the constant $C$ depends only on $\|\varphi_0\|_{L^\infty(M)}$, $n$, $\beta$, $\omega_{0}$, $F$ and $\overline{T}$. From $(\ref{207})$ and $(\ref{208})$, we prove the lemma.\QEDB

By Lemma \ref{205} and the fact that $\omega_\varepsilon>\gamma \omega_0$ for some uniform constant $\gamma$ (see inequality $(24)$ in \cite{CGP}), we have
\begin{equation}
e^{-h_C(t)}\omega_0 \leq\omega_{\varepsilon, j}\leq C_\varepsilon e^{h_C(t)}\omega_0
\end{equation}
on $(0,T]\times M$, where $C$ is a uniform constant and $C_\varepsilon$ depends on $\varepsilon$. We next prove the Calabi's $C^3$-estimates. Denote
\begin{eqnarray}S_{\varepsilon,j}=|\nabla_{\omega_0}\omega_{\varepsilon,j}(t)|^{2}_{\omega_{\varepsilon,j}(t)}=
g_{\varepsilon,j}^{i\bar{m}}g_{\varepsilon,j}^{k\bar{l}}g_{\varepsilon,j}^{p\bar{q}}\nabla_{0 i}(g_{\varepsilon,j}) _{k\bar{q}}\overline{\nabla}_{0 m}(g_{\varepsilon,j})_{ p\bar{l}}.
\end{eqnarray}

\begin{lem}\label{209} For any $0<T<\overline{T}_{\varepsilon,j}$ and $\varepsilon>0$, there exist constants $C_\varepsilon$ and $C$ such that for any $t\in(0,T]$, we have
\begin{eqnarray}S_{\varepsilon,j}\leq C_\varepsilon e^{h_C(t)},
\end{eqnarray}
where constant $C$ depends only on $\|\varphi_0\|_{L^\infty(M)}$, $n$, $\beta$, $\omega_{0}$, $F$ and $\overline{T}$, and constant $C_\varepsilon$ depends in addition on $\varepsilon$.
\end{lem}

{\bf Proof:}\ \ Straightforward calculations show that
\begin{equation}
\begin{split}
(\frac{\partial}{\partial t}-\Delta_{\varepsilon, j})S_{\varepsilon, j}(t)=&-|\nabla X|_{\varepsilon, j}^2-|\overline{\nabla} X|_{\varepsilon, j}^2+2Re \langle\nabla T-\nabla R_{m\omega_0}, X\rangle_{\omega_{\varepsilon, j}}\\
&+X^{i}_{pk}\overline{X^{j}_{ql}}(T^{p\bar{q}}g_{\varphi_\varepsilon i\bar{j}}g_{\varphi_{\varepsilon}}^{k\bar{l}}-g_{\varphi_{\varepsilon}}^{p\bar{q}}T_{i\bar{j}}g_{\varphi_{\varepsilon}}^{k\bar{l}}
+g_{\varphi_{\varepsilon}}^{p\bar{q}}g_{\varphi_{\varepsilon}i\bar{j}}T^{k\bar{l}}),
\end{split}
\end{equation}
where $X^{k}_{il}=(\nabla_{i}h\cdot h^{-1})^{k}_{\ l}$, $h^{i}_{\ k}=g_{0}^{i\overline{j}}g_{\varphi_{\varepsilon}\overline{j}k}$ and
\begin{equation}
\begin{split}
T_{i\bar{j}}=&(\frac{\partial g_{\varphi_{\varepsilon}}}{\partial t}+Ric(g_{\varphi_{\varepsilon}}))_{i\bar{j}}\\
=&Ric(g_0)_{i\bar{j}}+\log(\varepsilon^2+|s|_h^2)_{i\bar{j}}+F_{i\bar{j}}
+F''\varphi_i\varphi_{\bar{j}}+F'_i\varphi_{\bar{j}}+F'_{\bar{j}}\varphi_i+F'\varphi_{i\bar{j}}.
\end{split}
\end{equation}
Since
\begin{eqnarray*}
\nabla_m T_{i\bar{j}}&=&\nabla_{0 m} T_{i\bar{j}}-X_{ml}^s T_{s\bar{j}},\\
\nabla_{\varphi_{\varepsilon} p}R_{0\ l\overline{q}m}^{\ \beta}&=&\nabla_{0 p}R_{0\ l\overline{q}m}^{\ \beta}+X^{\beta}_{ps}R_{0\ l\overline{q}m}^{\ s}-X^{s}_{pl}R_{0\ s\overline{q}m}^{\ \beta}-X^{s}_{pm}R_{0\ l\overline{q}s}^{\ \beta},
\end{eqnarray*}
we have
\begin{equation}
\begin{split}
(\frac{\partial}{\partial t}-\Delta_{\varepsilon, j})S_{\varepsilon, j}(t)\leq& -|\nabla_{\varphi_{\varepsilon}} X|_{\varepsilon, j}^2-|\overline{\nabla}_{\varphi_{\varepsilon}} X|_{\varepsilon, j}^2+C_\varepsilon h_C(t)e^{h_C(t)}S_{\varepsilon,j}\\
&+C_\varepsilon h_C(t)e^{h_C(t)}+|\varphi_{ij}|_{\omega_{\varepsilon, j}}^2.
\end{split}
\end{equation}
Combining with the following evolution equations
\begin{equation}
\begin{split}
(\frac{\partial}{\partial t}-\Delta_{\varepsilon, j})\tr_{\omega_0}\omega_{\varepsilon, j}
\leq &-e^{-h_C(t)}S_{\varepsilon,j}+C_\varepsilon h_C(t)e^{h_C(t)},
\end{split}
\end{equation}
\begin{equation}
\begin{split}
(\frac{\partial}{\partial t}-\Delta_{\varepsilon, j})|\nabla \varphi_{\varepsilon, j}|_{\omega_0}^2
\leq & C_{\varepsilon} h_C(t)e^{h_C(t)}-e^{-h_C(t)}\sum_{j, p}\frac{|\varphi_{jp}|^2}{(1+\varphi_{p\bar{p}})(1+\varphi_{j\bar{j}})},
\end{split}
\end{equation}
we have
\begin{equation}
\begin{split}
&(\frac{\partial}{\partial t}-\Delta_{\varepsilon, j})(e^{-h_\alpha(t)}S_{\varepsilon,j}+ A_\varepsilon e^{-h_\beta(t)}\tr_{\omega_0}\omega_{\varepsilon, j}+e^{-h_\gamma(t)}|\nabla \varphi_{\varepsilon, j}|_{\omega_0}^2)\\
\leq &C h_C(t)e^{-h_C(t)}+e^{-h_C(t)}+C_\varepsilon e^{-\frac{h_\alpha(t)}{2}}S_{\varepsilon,j}-A_\varepsilon e^{-2h_\beta(t)}S_{\varepsilon,j}\\
&+e^{-h_\alpha(t)}|\varphi_{ij}|_{\omega_{\varepsilon, j}}^2- e^{-h_C(t)-h_\gamma(t)}|\varphi_{ij}|_{\omega_{\varepsilon, j}}^2.
\end{split}
\end{equation}
By choosing suitable $\alpha$, $\beta$ and $\gamma$ and using the maximum principle, we have
\begin{equation}
S_{\varepsilon,j}\leq C_\varepsilon e^{h_\alpha(t)},
\end{equation}
where $C_\varepsilon$ is independent of $j$ and $h_\alpha(t)$ is independent of $\varepsilon, j$.

By using the Schauder regularity theory and equation $(\ref{CMAE1})$, we get the high order estimates of $\varphi_{\varepsilon,j}(t)$.

\begin{pro}\label{210} For any $0<\eta<T<\overline{T}_{\varepsilon,j}$, $\varepsilon>0$ and $k\geq0$, there exists a constant $C_{k}$ depending only on $\eta$, $\overline{T}$, $\varepsilon$, $k$, $n$, $\beta$, $\omega_{0}$, $F$ and $\|\varphi_0\|_{L^\infty(M)}$, such that
\begin{eqnarray}\|\varphi_{\varepsilon,j}(t)\|_{C^{k}\big([\eta,T]\times M\big)}\leq C_{k}.
\end{eqnarray}
\end{pro}

\begin{pro}\label{21000000}
There exists a uniform $T$ such that the equation $(\ref{CMAE1})$ admits a unique solution on $[0,T]\times M$ for any $\varepsilon$ and $j$.
Furthermore,
\begin{eqnarray}\label{201}
\|\varphi_{\varepsilon,j}(t)-\varphi_{\varepsilon,l}(t)\|_{L^\infty([0,T]\times M)}\leq e^{\beta T}\|\varphi_{0,j}-\varphi_{0,l}\|_{L^\infty(M)}.
\end{eqnarray}
In particular, $\{\varphi_{\varepsilon,j}(t)\}$ satisfies
\begin{eqnarray}\label{2001}\lim\limits_{j,l\rightarrow\infty}\|\varphi_{\varepsilon,j}(t)-\varphi_{\varepsilon,l}(t)\|_{L^{\infty}([0,T]\times M)}=0.
\end{eqnarray}
\end{pro}

{\bf Proof.}\ \ Let $T_{\varepsilon,j}$ be the maximal existence interval of the equation $(\ref{CMAE1})$. We prove that $T_{\varepsilon,j}\geq \overline{T}$ for any $\varepsilon$ and $j$. If $T_{\varepsilon,j}<\overline{T}$, by Proposition $2.6$, we obtain $C^{\infty}$-estimates for $\varphi_{\varepsilon,j}(t)$ on $[0,T_{\varepsilon,j})$. Hence as $t\rightarrow T_{\varepsilon,j}$, $\varphi_{\varepsilon,j}(t)$ converge in $C^{\infty}$ to a smooth function $\varphi_{\varepsilon,j}(T_{\varepsilon,j})$ and thus we obtain a smooth solution to the equation $(\ref{CMAE1})$ on $[0,T_{\varepsilon,j}]$. But we can always find a smooth solution of the equation $(\ref{CMAE1})$ on some, possibly short, time interval with any initial K\"ahler potential. Applying this to $\varphi_{\varepsilon,j}(T_{\varepsilon,j})$, we obtain a solution of the equation $(\ref{CMAE1})$ on $[0, T_{\varepsilon,j}+\bar{\varepsilon})$ for $ \bar{\varepsilon}>0$. But this contradicts the definition of $T_{\varepsilon,j}$. Then $T_{\varepsilon,j}\geq\overline{T}$. Therefore, we conclude that there exists a uniform $T$ such that the equation $(\ref{CMAE1})$ admits a solution on $[0,T]\times M$ for any $\varepsilon$ and $j$. Next, we prove $(\ref{201})$.
\begin{equation}
\|\varphi_{\varepsilon, j}(t)\|_{L^\infty} \leq C,  \qquad for\ any\ t\in [0, T], \varepsilon >0\ and\ j\in \mathbb{N}^+.
\end{equation}
Set $\Psi_{\varepsilon, j, l}(t)= \varphi_{\varepsilon, j}(t)-\varphi_{\varepsilon, l}(t)$.
\begin{equation}
\begin{split}
\frac{\partial \Psi_{\varepsilon, j, l}(t)}{\partial t}=&\log\frac{(\omega_0+\sqrt{-1}\partial\bar{\partial}\varphi_{\varepsilon,l}(t)+\sqrt{-1}\partial\bar{\partial}\Psi_{\varepsilon, j, l}(t))^{n}}{(\omega_0+\sqrt{-1}\partial\bar{\partial}\varphi_{\varepsilon,l}(t))^{n}}\\
&+F(\varphi_{\varepsilon, j}(t), z)-F(\varphi_{\varepsilon, l}(t), z).\\
\end{split}
\end{equation}
Since \begin{equation}
|F(\varphi_{\varepsilon, j}(t), z)-F(\varphi_{\varepsilon, l}(t), z)|\leq K|\varphi_{\varepsilon, j}(t)-\varphi_{\varepsilon, l}(t)|,
\end{equation}
where $K$ only depends on $\|\varphi_0\|_{L^\infty}, \omega_0, \beta, n, F$ and $T$,
\begin{equation}
\frac{\partial \Psi_{\varepsilon, j, l}(t)}{\partial t}\leq\log\frac{(\omega_{\varepsilon, l}+\sqrt{-1}\partial\bar{\partial}\Psi_{\varepsilon, j, l}(t))^{n}}{\omega_{\varepsilon, l}^n}+K|\Psi_{\varepsilon, j, l}(t)|.\\
\end{equation}
Let $\widetilde{\Psi}_{\varepsilon, j, l}(t)= e^{-Kt}\Psi_{\varepsilon, j, l}(t)-\delta t$. Assume that the maximum point of $\widetilde{\Psi}_{\varepsilon, j, l}(t)$ is $(t_0, x_0)$. If $t_0 =0$, then
\begin{equation}
\Psi_{\varepsilon, j, l}(t)\leq e^{KT}\sup_M\Psi_{\varepsilon, j, l}(0)+\delta Te^{KT}.
\end{equation}
If $t_0> 0$ and $\Psi_{\varepsilon, j, l}(t_0, x_0)\leq 0$, we have
\begin{equation}
\Psi_{\varepsilon, j, l}(t)\leq\delta T e^{KT}.
\end{equation}
If $t_0> 0$ and $\Psi_{\varepsilon, j, l}(t_0, x_0)> 0$, then
\begin{equation}
\begin{split}
&\frac{\partial}{\partial t}(e^{-Kt}\Psi_{\varepsilon, j, l}(t)-\delta t)\\
\leq& -Ke^{-Kt}\Psi_{\varepsilon, j, l}(t)+e^{-Kt}\log\frac{(e^{-Kt}\omega_{\varepsilon, l}+\sqrt{-1}\partial\bar{\partial}(e^{-Kt}\Psi_{\varepsilon, j, l}(t)-\delta t))^{n}}{(e^{-Kt}\omega_{\varepsilon, l})^n}\\
&+Ke^{-Kt}|\Psi_{\varepsilon, j, l}(t)|-\delta.
\end{split}
\end{equation}
By the maximum principle, at $(t_0, x_0)$, we have $0\leq \frac{\partial \widetilde{\Psi}_{\varepsilon, j, l}}{\partial t}(t_0, x_0)\leq -\delta$, which is impossible. Hence
\begin{equation}\label{2016080801}
\Psi_{\varepsilon, j, l}(t)\leq e^{KT}\|\varphi_{0, j}-\varphi_{0, l}\|_{L^\infty}+\delta Te^{KT}.
\end{equation}
By the same arguments, we can get the lower bound of $\Psi_{\varepsilon,j,l}(t)$,
\begin{equation}\label{2016080802}
\Psi_{\varepsilon, j, l}(t)\geq -e^{KT}\|\varphi_{0, j}-\varphi_{0, l}\|_{L^\infty(M)}-\delta Te^{KT}.
\end{equation}
Combining $(\ref{2016080801})$ and $(\ref{2016080802})$, we have
\begin{equation}
\|\varphi_{\varepsilon, j}(t)-\varphi_{\varepsilon, l}(t)\|_{L^\infty([0, T]\times M)}\leq e^{KT}\|\varphi_{0, j}(t)-\varphi_{0, l}(t)\|_{L^\infty(M)}+\delta Te^{KT}.
\end{equation}
Let $\delta \rightarrow 0$, we have
\begin{equation}
\|\varphi_{\varepsilon, j}(t)-\varphi_{\varepsilon, l}(t)\|_{L^\infty([0, T]\times M)}\leq e^{KT}\|\varphi_{0, j}(t)-\varphi_{0, l}(t)\|_{L^\infty(M)}.
\end{equation}
By the similar arguments as above, we know that the solution to equation $(\ref{CMAE2})$ must be unique. \QEDB

By $(\ref{2001})$, $\varphi_{\varepsilon,j}(t)$ converges to $\varphi_\varepsilon(t)\in L^\infty([0,T]\times M)$ uniformly in $L^\infty([0,T]\times M)$-topology. For any $0<\delta<T$ and $\varepsilon>0$, $\varphi_{\varepsilon,j}(t)$ is uniformly bounded (depends on $\varepsilon$) in $C^{\infty}([\delta,T]\times M)$. Therefore $\varphi_{\varepsilon,j}(t)$ converges to $\varphi_\varepsilon(t)$ in $C^\infty([\delta,T]\times M)$-topology. Hence for any $\varepsilon>0$, $\varphi_\varepsilon(t)\in C^\infty((0,T]\times M)$.

\begin{pro}\label{211} For any $\varepsilon>0$, $\varphi_\varepsilon(t)\in C^0([0,T]\times M)$ and
\begin{eqnarray}\label{212}\lim\limits_{t\rightarrow0^+}\|\varphi_{\varepsilon}(t)-\varphi_{0}\|_{L^{\infty}(M)}=0.
\end{eqnarray}
\end{pro}

{\bf Proof:}\ \ For any $(t,z)\in(0,T]\times M$,
\begin{eqnarray}\nonumber|\varphi_{\varepsilon}(t,z)-\varphi_{0}(z)|&\leq&|\varphi_{\varepsilon}(t,z)-\varphi_{\varepsilon,j}(t,z)|
+|\varphi_{\varepsilon,j}(t,z)-\varphi_{0,j}(z)|\\
&\ &+|\varphi_{0,j}(z)-\varphi_{0}(z)|.
\end{eqnarray}
Since $\varphi_{\varepsilon,j}(t)$ is a Cauchy sequence in $L^\infty([0,T]\times M)$,
\begin{eqnarray}\lim\limits_{j\rightarrow\infty}\|\varphi_{\varepsilon}(t,z)-\varphi_{\varepsilon,j}(t,z)\|_{L^{\infty}([0,T]\times M)}=0.
\end{eqnarray}
From $(\ref{00000})$, we have
\begin{eqnarray}
\lim\limits_{j\rightarrow\infty}\|\varphi_{0,j}(z)-\varphi_{0}(z)\|_{L^{\infty}(M)}=0.
\end{eqnarray}
For any $\epsilon>0$, there exists an $N$ such that for any $j>N$,
\begin{eqnarray*}\sup\limits_{[0,T]\times M}|\varphi_{\varepsilon}(t,z)-\varphi_{\varepsilon,j}(t,z)|&<&\frac{\epsilon}{3},\\
\sup\limits_{M}|\varphi_{0,j}(z)-\varphi_{0}(z)|&<&\frac{\epsilon}{3}.
\end{eqnarray*}
On the other hand, fixed such $j$, there exists $0<\delta<T$ such that
\begin{eqnarray}\sup\limits_{[0,\delta]\times M}|\varphi_{\varepsilon,j}(t,z)-\varphi_{0,j}|<\frac{\epsilon}{3}.
\end{eqnarray}
Combining the above estimates together, for any $t\in[0,\delta]$ and $z\in M$,
\begin{eqnarray}|\varphi_{\varepsilon}(t,z)-\varphi_{0}(z)|<\epsilon.
\end{eqnarray}
This completes the proof of the lemma.\QEDB

\begin{pro} \label{214}$\varphi_\varepsilon(t)$ is the unique solution to the parabolic Monge-Amp\`ere equation $(\ref{3})$ in the space of $C^0\big([0,T]\times M\big)\cap C^\infty\big((0,T]\times M\big)$.
\end{pro}

\medskip

{\bf Proof:}\ \ By proposition \ref{211}, we only need to prove the uniqueness. Suppose there exists another solution $\tilde{\varphi}_\varepsilon(t)\in C^0\big([0,T]\times M\big)\cap C^\infty\big((0,T]\times M\big)$ to the Monge-Amp\`ere equation $(\ref{3})$. Let $\psi_\varepsilon(t)=\tilde{\varphi}_\varepsilon(t)-\varphi_\varepsilon(t)$. Then
\begin{eqnarray*}
\begin{cases}
  \frac{\partial \psi_{\varepsilon}(t)}{\partial t}=\log\frac{\big(\omega_{0}+\sqrt{-1}\partial\bar{\partial}\varphi_{\varepsilon}(t)+\sqrt{-1}\partial\bar{\partial}\psi_{\varepsilon}(t)\big)^{n}}
  {(\omega_{0}+\sqrt{-1}\partial\bar{\partial}\varphi_{\varepsilon}(t))^{n}}+F(\tilde{\varphi}_{\varepsilon}(t), z)-F(\varphi_{\varepsilon}(t), z),\\
  \\
  \psi_{\varepsilon}(0)=0.\\
  \end{cases}
\end{eqnarray*}
By the same arguments as that in the proof of Proposition \ref{21000000},
\begin{eqnarray*}
\|\psi_{\varepsilon}(t)\|_{L^\infty([0,T]\times M)}\leq e^{\beta T}\|\psi_{\varepsilon}(0)\|_{L^\infty(M)}=0.
\end{eqnarray*}
Hence $\psi_{\varepsilon}(t)=0$, that is, $\tilde{\varphi}_\varepsilon(t)=\varphi_\varepsilon(t)$.\QEDB

\section{The existence of the conical parabolic complex Monge-Amp\`ere equation with weak initial data}
\setcounter{equation}{0}

In this section, we study the existence of the conical parabolic complex Monge-Amp\`ere equation $(\ref{1})$ by the smooth approximation of equations $(\ref{3})$. We also prove the uniqueness and regularity of the equation $(\ref{1})$.

By Lemma \ref{2016080701}, Lemma \ref{205}, Proposition \ref{210} and Proposition \ref{21000000}, we conclude that there exist constants $C_1$ and $C_{2}$ depending only on $\|\varphi_0\|_{L^\infty(M)}$, $\beta$, $n$, $\omega_{0}$, $F$ and $T$, such that for any $\varepsilon>0$,
\begin{eqnarray}\label{219}&\ &\|\varphi_{\varepsilon}(t)\|_{L^\infty([0,T]\times M)}\leq C_1,\\
\label{220}&\ &e^{-h_{C_2}(t)}\omega_\varepsilon\leq\omega_{\varepsilon}(t)\leq e^{h_{C_2}(t)}\omega_\varepsilon\ \ on\ (0,T]\times M.
\end{eqnarray}

We first prove the local uniform Calabi's $C^3$-estimates and curvature estimates along the equation $(\ref{CMAE1})$. Our proofs are similar as that in Section $2$ of \cite{JWLXZ1}, but we need some arguments to handle the terms from $F$.
\begin{lem} For any $0<\eta<T$ and $B_r(p)\subset\subset (M\setminus D)$, there exist constants $C'$ and $C''$ such that for any $\varepsilon>0$ and $j\in\mathbb{N}^+$,
\begin{eqnarray*}S_{\varepsilon,j}&\leq&\frac{C'}{r^{2}},\\
|Rm_{\varepsilon,j}|_{\omega_{\varepsilon,j}(t)}^{2}&\leq&\frac{C''}{r^{4}}
\end{eqnarray*}
on $[\eta,T]\times B_{\frac{r}{2}}(p)$, where constants $C'$ and $C''$ depend only on $\|\varphi_0\|_{L^\infty(M)}$, $n$, $\beta$, $\eta$, $T$, $\omega_0$, $F$ and $dist_{\omega_0}(B_r(p),D)$.\end{lem}

{\bf Proof:}\ \ By Lemma \ref{205}, there exists a uniform constat $C$ depending only on $\|\varphi_0\|_{L^\infty(M)}$, $n$, $\beta$, $T$, $\omega_0$, $F$ and $dist_{\omega_0}(B_r(p),D)$, such that
\begin{eqnarray}e^{-h_{C}(t)}\omega_0\leq\omega_{\varepsilon,j}(t)\leq e^{h_{C}(t)}\omega_0,\ \ on\ \ B_r(p)\times (0,T].
\end{eqnarray}

Let $r=r_0>r'_1>\frac{r}{2}$ and $\psi$ be a nonnegative $C^\infty$ cut-off function that is identically equal to $1$ on $\overline{B_{r_{1}(p)}}$ and vanishes outside $B_{r'}(p)$. We may assume that
\begin{eqnarray}|\partial\psi|_{\omega_{0}}^{2}\leq\frac{C}{r^{2}}\ \ \ and\ \ \ |\sqrt{-1}\partial\bar{\partial}\psi|_{\omega_{0}}\leq\frac{C}{r^{2}}.\end{eqnarray}
Straightforward calculations show that
\begin{eqnarray}(\frac{\partial}{\partial t}-\Delta_{\varepsilon,j})(\psi^2e^{-h_{\alpha}(t)}S_{\varepsilon,j})\leq \frac{C}{r^2} e^{-\frac{h_{\alpha}(t)}{2}}S_{\varepsilon,j}+e^{-h_{\alpha}(t)}|\varphi_{ij}|_{\omega_{\varepsilon,j}}+C,
\end{eqnarray}
\begin{equation}
(\frac{\partial}{\partial t}-\Delta_{\varepsilon, j})(A e^{-h_\beta(t)}\tr_{\omega_0}\omega_{\varepsilon, j})
\leq -A e^{-h_c(t)-h_\beta(t)}S_{\varepsilon,j}+AC e^{\frac{-h_\beta(t)}{2}}+A,
\end{equation}
\begin{equation}
(\frac{\partial}{\partial t}-\Delta_{\varepsilon, j})(e^{-h_\gamma(t)}|\nabla \varphi_{\varepsilon, j}|_{\omega_0}^2)
\leq -e^{-h_c(t)-h_\gamma(t)}|\varphi_{ij}|_{\omega_{\varepsilon, j}}^2+e^{-\frac{h_\gamma(t)}{2}}+C.
\end{equation}
By choosing sufficiently large $\alpha$, $\beta$, $\gamma$, $A$ and using the maximum principle, we conclude that
\begin{eqnarray*}S_{\varepsilon,j}&\leq&\frac{C'}{r^{2}}e^{h_{\alpha}(t)}\ \ \ on\ \ (0,T]\times B_{r'_{1}}(p).
\end{eqnarray*}
Hence $\|\varphi_{\varepsilon,j}(t)\|_{C^{2,\alpha}}$ is uniformly bounded on $[\eta_{1}, T]\times\overline{ B_{r_{1}}}$ with $0<\eta_{1}<\eta$ and $\frac{r}{2}<r_{1}<r'_{1}$.

Now we prove that $|Rm_{\varepsilon,j}|^{2}_{\omega_{\varepsilon,j}(t)}$ is uniformly bounded. On $[\eta_1, T]\times\overline{B_{r_1}}$, since $\|\varphi_{\varepsilon,j}(t)\|_{C^{2,\alpha}}$ is uniformly bounded,
\begin{equation}
\begin{split}
&(\frac{\partial}{\partial t}-\Delta_{\varepsilon, j})|Rm_{\varepsilon,j}|_{\omega_{\varepsilon, j}}^2\\
\leq&C|Rm_{\varepsilon,j}|^3+Ch_C(t)e^{h_C(t)}|Rm_{\varepsilon,j}|_{\omega_{\varepsilon,j}}^2+\frac{Ch_C(t)e^{h_C(t)}}{r^2}|Rm_{\varepsilon,j}|_{\omega_{\varepsilon,j}}\\
&+\frac{Ch_C(t)e^{h_C(t)}}{r^2}|\varphi_{ik}|_{\omega_0}|Rm_{\varepsilon,j}|_{\omega_{\varepsilon,j}}+Ce^{h_C(t)}|\varphi_{ik}|_{\omega_0}^2|Rm_{\varepsilon,j}|_{\omega_{\varepsilon,j}}\\
\leq & C(|Rm_{\varepsilon,j}|_{\omega_{\varepsilon,j}}^3+1+\frac{1}{r^2}|Rm_{\varepsilon,j}|_{\omega_{\varepsilon,j}})-|\nabla Rm_{\varepsilon,j}|_{\omega_{\varepsilon,j}}^2-|\overline{\nabla}Rm_{\varepsilon,j}|_{\omega_{\varepsilon,j}}^2,
\end{split}
\end{equation}
where $C$ depends on $\delta$ and $T$.

We fix a smaller radius $r_{2}$ satisfying $r_{1}>r_{2}>\frac{r}{2}$. Let $\rho$ be a cut-off function identically equal to $1$ on $\overline{B_{r_{2}}}(p)$ and identically equal to $0$ outside $B_{r_{1}}$. We also let $\rho$ satisfy
$$|\partial\rho|_{\omega_{0}}^{2},\ |\sqrt{-1}\partial\bar{\partial}\rho|_{\omega_{0}}\leq\frac{C}{r^{2}}$$
for some uniform constant $C$. From the former part we know that $S_{\varepsilon,j}$ is bounded on $[\eta_{1},T]\times B_{r_{1}}(p)$. Let $K=\frac{\hat{C}}{r^{2}}$, $\hat{C}$ be constants which are large enough such that $\frac{K}{2}\leq K-S_{\varepsilon,j}\leq K$. We consider
\begin{eqnarray}\label{3.20.11}F_{\varepsilon,j}=\rho^{2}e^{-\frac{2\alpha}{t-\eta_{1}}}\frac{|Rm_{\varepsilon,j}|^{2}_{\omega_{\varepsilon,j}(t)}}{K-S_{\varepsilon,j}}
+Ae^{-\frac{2\beta}{t-\eta_{1}}}S_{\varepsilon,j}.\end{eqnarray}
We only consider an inner point $(t_{0},x_{0})$ which is a maximum point of $F_{\varepsilon,j}$ achieved on $[0,T]\times\overline{B_{r_{1}}(p)}$. By the similar arguments as that in Lemma $3.1$ of \cite{JWLXZ1}, we have
\begin{equation}
\begin{split}
(\frac{d}{dt}-\Delta_{\varepsilon, j})F_{\varepsilon, j}\leq &-\frac{A}{2}e^{-\frac{2\beta}{t-\delta_1}}Q+Ce^{-\frac{\alpha}{t-\delta_1}}Q+Ae^{-\frac{2\beta}{t-\delta_1}}\frac{C}{r^2},
\end{split}
\end{equation}
where $Q=|\nabla X|^{2}+|\bar{\nabla} X|^{2}$. Now we choose $\alpha=2\beta$ and $\frac{A}{2}=C+1$, then
$$e^{-\frac{\alpha}{t-\eta_1}}Q\leq \frac{C}{r^2}$$
at $(t_{0},x_{0})$. This implies that
\begin{eqnarray}\label{3.20.21}|Rm_{\varepsilon,j}|_{\omega_{\varepsilon,j}(t)}^{2}\leq\frac{C_{1}}{r^{4}}e^{\frac{C_{2}}{t-\eta_{1}}},\end{eqnarray}
where $C_{1}$ and $C_{2}$ depend only on $\|\varphi_0\|_{L^\infty(M)}$, $n$, $\beta$, $T$, $dist_{\omega_{0}}(B_r(p),D)$, $\eta_{1}$, $F$ and $\omega_0$. Hence we prove the lemma.\QEDB

\medskip

Using the standard parabolic Schauder regularity theory \cite{GLIE}, we obtain the following proposition.

\begin{pro}\label{217} For any $0<\eta<T$, $k\in\mathbb{N}^{+}$ and $B_r(p)\subset\subset (M\setminus D)$, there exists a constant $C_{\eta,T,k,p,r}$ which depends only on $\|\varphi_0\|_{L^\infty(M)}$, $n$, $\beta$, $\eta$, $k$, $T$, $dist_{\omega_{0}}(B_r(p),D)$, $F$ and $\omega_0$, such that for any $\varepsilon>0$ and $j\in\mathbb{N}^+$,
\begin{eqnarray}\|\varphi_{\varepsilon,j}(t)\|_{C^{k}\big([\eta,T]\times B_r(p)\big)}\leq C_{\eta,T,k,p,r}.
\end{eqnarray}
\end{pro}

Through a further observation to equation $(\ref{3})$, we prove the monotonicity of $\varphi_\varepsilon(t)$ with respect to $\varepsilon$.

\begin{pro}\label{218} For any $(t,x)\in [0,T]\times M$, $\varphi_\varepsilon(t,x)$ is monotonously  decreasing as $\varepsilon\searrow0$.
\end{pro}

{\bf Proof:}\ \ For any $\varepsilon_1<\varepsilon_2$, let $\psi_{1,2}(t)=\varphi_{\varepsilon_1}(t)-\varphi_{\varepsilon_2}(t)$. Then we have
\begin{eqnarray}\nonumber&\ &\frac{\partial \psi_{1,2}}{\partial t}\leq \log\frac{(\omega_{\varphi_{\varepsilon_2}}+\sqrt{-1}\partial\bar{\partial}\psi_{1,2}(t))^{n}}
{\omega_{\varphi_{\varepsilon_2}}^{n}}+K|\psi_{1, 2}|.
\end{eqnarray}

By the similar arguments as that in Proposition \ref{21000000}, for any $(t,x)\in [0,T]\times M$ and $\delta>0$, we have
\begin{eqnarray}\psi_{1,2}(t,x)\leq e^{K t}\sup\limits_{M}\psi_{1,2}(0,x)+Te^{K T}\delta=Te^{K T}\delta.
\end{eqnarray}
Let $\delta\rightarrow0$, we conclude that $\varphi_{\varepsilon_1}(t,x)\leq\varphi_{\varepsilon_2}(t,x)$.\QEDB

For any $[\delta, T]\times K\subset\subset (0,T]\times (M\setminus D)$ and $k\geq0$, $\|\varphi_{\varepsilon,j}(t)\|_{C^{k}([\delta,T]\times K)}$ is uniformly bounded by Proposition \ref{217}. Let $j$ approximate to $\infty$, we obtain that $\|\varphi_{\varepsilon}(t)\|_{C^{k}([\delta,T]\times K)}$ is uniformly bounded. Then let $\delta$ approximate to $0$ and $K$ approximate to $M\setminus D$, by the diagonal rule, we get a sequence $\{\varepsilon_i\}$, such that $\varphi_{\varepsilon_i}(t)$ converges in $C^\infty_{loc}$ topology on $(0,T]\times (M\setminus D)$ to a function $\varphi(t)$ that is smooth on $C^\infty\big((0,T]\times(M\setminus D)\big)$ and satisfies the equation
\begin{eqnarray}\label{223}  \frac{\partial \varphi(t)}{\partial t}=\log\frac{(\omega_{0}+\sqrt{-1}\partial\bar{\partial}\varphi(t))^{n}}{\omega_{0}^{n}}+F_{0}+\beta\varphi(t)
  +\log|s|_{h}^{2(1-\beta)}
\end{eqnarray}
on $(0,T]\times (M\setminus D)$. Since $\varphi_\varepsilon(t)$ is decreasing as $\varepsilon\rightarrow0$, we conclude that $\varphi_{\varepsilon}(t)$ converges in $C^\infty_{loc}$ topology on $(0,T]\times (M\setminus D)$ to $\varphi(t)$. Combining the above arguments with
$(\ref{219})$ and $(\ref{220})$, for any $T>0$, we have
\begin{eqnarray}\label{221}&\ &\|\varphi(t)\|_{L^\infty\big((0,T]\times (M\setminus D)\big)}\leq C_1,\\
\label{222}&\ &e^{-h_{C_2}(t)}\omega_\beta\leq\omega(t)\leq e^{h_{C_2}(t)}\omega_\beta\ \ on\ (0,T]\times (M\setminus D),
\end{eqnarray}
where $\omega(t)=\omega_0+\sqrt{-1}\partial\overline{\partial}\varphi(t)$, constants $C_1$ and $C_{2}$ depend only on $\|\varphi_0\|_{L^\infty(M)}$, $\beta$, $n$, $\omega_{0}$, $F$ and $T$. By the similar arguments as that in \cite{JWLXZ}, we obtain the following proposition.

\begin{pro}\label{2188} For any $t\in(0,T]$, $\varphi(t)$ is H\"older continuous on $M$ with respect to the metric $\omega_0$.
\end{pro}

Next, by using the monotonicity of $\varphi_{\varepsilon}(t)$ with respect to $\varepsilon$ and constructing a auxiliary function, we prove the continuity of $\varphi(t)$ as $t\rightarrow0^{+}$.

\begin{pro}\label{101} $\varphi(t)\in C^0([0,T]\times M)$ and
\begin{eqnarray}\label{102}\lim\limits_{t\rightarrow0^+}\|\varphi(t)-\varphi_{0}\|_{L^{\infty}(M)}=0.
\end{eqnarray}
\end{pro}

{\bf Proof:}\ \ Through the above arguments, we only need prove the limit $(\ref{102})$. By the monotonicity of $\varphi_{\varepsilon}(t)$ with respect to $\varepsilon$, for any $\epsilon>0$ and $(t,z)\in(0,T]\times M$,
\begin{eqnarray}\label{2210}\nonumber\varphi(t,z)-\varphi_0(z)&\leq&\varphi_{\varepsilon_{1}}(t,z)-\varphi_0(z)\\\nonumber
&\leq&|\varphi_{\varepsilon_1}(t,z)-\varphi_{\varepsilon_1,j}(t,z)|
+|\varphi_{\varepsilon_1,j}(t,z)-\varphi_{0,j}(z)|\\
&\ &+|\varphi_{0,j}(z)-\varphi_{0}(z)|.
\end{eqnarray}
By the same arguments as that in Proposition \ref{211}, for any $t\in(0,\delta_1]$ and $z\in M$,
\begin{eqnarray}\label{201808}\varphi(t,z)-\varphi_{0}(z)<\epsilon.
\end{eqnarray}
On the other hand, by S. Kolodziej's results \cite{K}, there exists a smooth solution $u_{\varepsilon,j}$ to the equation
\begin{eqnarray}(\omega_0+\sqrt{-1}\partial\bar{\partial}u_{\varepsilon, j})^n=e^{-F(0, z)-K\varphi_{0, j}+\widehat{C}}\frac{\omega_0^n}{(\varepsilon^2+|s|_h^2)^{1-\beta}},
\end{eqnarray}
and $u_{\varepsilon,j}$ satisfies \begin{eqnarray}\|u_{\varepsilon,j}\|_{L^\infty(M)}\leq C,
\end{eqnarray}
where $\hat{C}$ is a uniform normalization constant independent of $\varepsilon$ and $j$, the constant $C$ depends only on $\|\varphi_0\|_{L^\infty(M)}$, $\beta$ and $F$.

We define
\begin{eqnarray}\Psi_{\varepsilon, j}(t)=(1-te^{Kt})\varphi_{0, j}+te^{Kt}u_{\varepsilon, j}+ h(t)e^{Kt}-m\|\varphi_{0, j}\|t-l|\widehat{C}|t,
\end{eqnarray}
where
\begin{equation}
\begin{split}
h(t)=&-t\|\varphi_{0, j}\|_{L^\infty}-t\|u_{\varepsilon, j}\|_{L^\infty}+n(t\log t-t)e^{-Kt}\\
&+Kn\int_0^te^{-Ks}s\log sds-\frac{\widehat{C}}{K}e^{-Kt}+\frac{\widehat{C}}{K},
\end{split}
\end{equation}
$m,l$ should be choosen later. Straightforward calculations show that
\begin{equation}
\begin{split}
\frac{\partial \Psi_{\varepsilon, j}}{\partial t}-F(\Psi_{\varepsilon, j}(t), z)\leq &K\|\varphi_{0,j}\|_{L^{\infty}(M)}+Kmt\|\varphi_{0, j}\|_{L^{\infty}(M)}+Klt|\widehat{C}|-\widehat{C}+n\log t\\
&+2\widehat{C}e^{Kt}-F(0,z)-m\|\varphi_{0, j}\|_{L^{\infty}(M)}-l|\widehat{C}|+nKt\\
\leq &-F(0,z)+n\log t+\widehat{C}+nKt-K\varphi_{0, j},
\end{split}
\end{equation}
when $t$ is sufficiently small and $m$, $l$ are sufficiently large.
Combining the above inequalities,
\begin{equation}
\begin{split}
(\omega_0+\sqrt{-1}\partial\bar{\partial}\Psi_{\varepsilon, j}(t))^{n}\geq &t^ne^{nKt}(\omega_0+\sqrt{-1}\partial\bar{\partial}u_{\varepsilon, j}(t))^{n}\\
\geq &e^{\frac{\partial \Psi_{\varepsilon, j}}{\partial t}-F(\Psi_{\varepsilon, j}(t), z)}\frac{\omega_0^n}{(\varepsilon^2+|s|_h^2)^{1-\beta}}.
\end{split}
\end{equation}
This equation is equivalent to
\begin{equation}
\left \{\begin{split}
\frac{\partial \Psi_{\varepsilon, j}}{\partial t}\leq & \log\frac{(\omega_0+\sqrt{-1}\partial\bar{\partial}\Psi_{\varepsilon, j}(t))^{n}}{\omega_{0}^{n}}+F(\Psi_{\varepsilon, j}(t), z)+\log(\varepsilon^2+|s|_h^2)^{1-\beta},\\
\Psi(0)=&\varphi_{0, j}.\\
\end{split}
\right.
\end{equation}
By the similar arguments as that in the proof of Proposition \ref{21000000}, for any $(t,z)\in[0,T]\times M$,
\begin{eqnarray}\varphi_{\varepsilon, j}(t)-\varphi_{0, j}\geq -te^{Kt}\varphi_{0, j}+te^{Kt}u_{\varepsilon, j}-mt\|\varphi_{0, j}\|-lt|\widehat{C}|+h(t)e^{Kt}.
\end{eqnarray}
Let $j\rightarrow\infty$ and then $\varepsilon\rightarrow0$, we have
\begin{eqnarray}\varphi(t,z)-\varphi_0(z)\geq-Cte^{K t}-Ct+h(t)e^{K t}.
\end{eqnarray}
There exists $\delta_2$ such that for any $t\in[0,\delta_2]$,
\begin{eqnarray}\varphi(t,z)-\varphi_0(z)>-\epsilon.
\end{eqnarray}
Let $\delta=\min(\delta_1,\delta_2)$, then for any $t\in(0,\delta]$ and $z\in M$,
\begin{eqnarray}-\epsilon<\varphi(t,z)-\varphi_{0}(z)<\epsilon.
\end{eqnarray}
This completes the proof of the proposition.\QEDB

Now we are ready to prove the uniqueness of the conical parabolic complex Monge-Amp\`ere equation $(\ref{1})$ starting with $\varphi_0\in\mathcal{E}_{p}(M,\omega_{0})$ for some $p>1$.

\begin{thm}\label{228} Let $\varphi_i(t)\in C^0\big([0,T]\times M\big)\bigcap C^\infty\big((0,T]\times(M\setminus D)\big)$ $(i=1,2)$ be two solutions to the parabolic Monge-Amp\`ere equation
\begin{eqnarray}
  \frac{\partial \varphi_i(t)}{\partial t}=\log\frac{(\omega_{0}+\sqrt{-1}\partial\bar{\partial}\varphi_i(t))^{n}}{\omega_{0}^{n}}+F(\varphi_i(t),z)
  +\log|s|_{h}^{2(1-\beta)}
\end{eqnarray}
on $(0,T]\times (M\setminus D)$. If $\varphi_i$ $(i=1,2)$ satisfy
\begin{itemize}
  \item For any $0<\delta<T$, there exists a uniform constant $C$ such that
\begin{eqnarray*}C^{-1}\omega_\beta\leq \omega_{0}+\sqrt{-1}\partial\bar{\partial}\varphi_i(t)\leq C\omega_\beta
\end{eqnarray*}
on $[\delta,T]\times (M\setminus D)$;
  \item On $[\delta, T]$, there exist constants $\alpha>0$ and $C^{*}$ such that $\varphi_i(t)$ is $C^{\alpha}$ on $M$ with respect to $\omega_{0}$ and $\| \frac{\partial\varphi_i(t)}{\partial t}\|_{L^{\infty}(M\setminus D)}\leqslant C^{*}$;
  \item $\lim\limits_{t\rightarrow0^{+}}\|\varphi_i(t)-\varphi_{0}\|_{L^{\infty}(M)}=0$.
  \end{itemize}
Then $\varphi_1=\varphi_2$.
\end{thm}

{\bf Proof:}\ \ We apply Jeffres' trick \cite{TJEF} in the parabolic case. For any $0<t_1<T$ and $a>0$, let $\phi_1(t)=\varphi_1(t)+a|s|_h^{2q}$, where $0<q<1$ is determined later. The evolution equation of $\phi_1$ is
\begin{eqnarray*}
  \frac{\partial \phi_1(t)}{\partial t}=\log\frac{(\omega_{0}+\sqrt{-1}\partial\bar{\partial}\varphi_1(t))^{n}}{\omega_{0}^{n}}+F(\varphi_1(t), z)+\log|s|_{h}^{2(1-\beta)}.
\end{eqnarray*}
Denote $\psi(t)=\phi_1(t)-\varphi_2(t)$ and $\hat{\Delta}=\int_0^1 g_{s\varphi_1+(1-s)\varphi_2}^{i\bar{j}}\frac{\partial^2}{\partial z^i\partial\bar{z}^j}ds$, $\psi(t)$ evolves along the following equation
\begin{eqnarray*}
  \frac{\partial \psi(t)}{\partial t}=\hat{\Delta}\psi(t)-a\hat{\Delta}|s|_h^{2q}+F(\varphi_1(t), z)-F(\varphi_2(t), z).\end{eqnarray*}
By the arguments as that in \cite{JWLXZ}, we obtain
\begin{eqnarray*}
  \frac{\partial \psi(t)}{\partial t}\leq\hat{\Delta}\psi(t)+K|\psi|+aC.
\end{eqnarray*}

Let $\tilde{\psi}=e^{-K (t-t_1)}\psi+\frac{aC}{K}e^{-K (t-t_1)}-\epsilon (t-t_1)$. By choosing suitable $0<q<1$, we can  assume that the space maximum of $\tilde{\psi}$ on $[t_1,T]\times M$ is attained away from $D$. Let $(t_0,x_0)$ be the maximum point. If $t_0>t_1$ and $\psi(t_{0},x_{0})>0$,  by the maximum principle, at $(t_0,x_0)$, we have
\begin{eqnarray*}
 0\leq (\frac{\partial }{\partial t}-\hat{\Delta})\tilde{\psi}(t)\leq -\epsilon,
\end{eqnarray*}
which is impossible, hence $t_0=t_1$ or $t_0>t_1$ and $\psi(t_{0},x_{0})\leq0$. Then for $(t,x)\in [t_1,T]\times M$, we obtain
\begin{eqnarray*}
\psi(t,x)&\leq& e^{K T}\|\varphi_1(t_1,x)-\varphi_2(t_1,x)\|_{L^\infty(M)}\\
&\ &+aCe^{K T}+\epsilon Te^{K T}.
\end{eqnarray*}
Let $a\rightarrow0$ and then $t_1\rightarrow0^+$, we get
\begin{eqnarray*}
\varphi_1(t)-\varphi_2(t)\leq \epsilon Te^{\beta T}.
\end{eqnarray*}
It shows that $\varphi_1(t)\leq\varphi_2(t)$ after we let $\epsilon\rightarrow0$. By the same reason we have $\varphi_2(t)\leq\varphi_1(t)$, then we prove that $\varphi_1(t)=\varphi_2(t)$.\QEDB

Next, we prove the uniform gradient estimates and Laplacian estimates of $\dot{\phi}_{\varepsilon,j}(t)$ on $M\times[\eta,T]$.
\begin{lem}\label{1.8.6}  For any $0<\eta<T$, there exists a uniform constant $C$ depending only on $\|\varphi_0\|_{L^\infty(M)}$, $n$, $\beta$, $\eta$, $\omega_{0}$, $F$ and $T$, such that for any $t\in[\eta,T]$, $\varepsilon>0$ and $j\in\mathbb{N}^+$, we have
\begin{eqnarray}\label{1.5.2}|\nabla \dot{\varphi}_{\varepsilon,j}(t)|_{\omega_{\varepsilon,j}(t)}^2&\leq& C,\\
\label{1.5.3}|\Delta_{\varepsilon,j} \dot{\varphi}_{\varepsilon,j}(t)|&\leq& C.\end{eqnarray}
\end{lem}

{\bf Proof:}\ \ As the computations in \cite{JWLXZ}, we have
\begin{equation}
(\frac{d}{d t}-\Delta_{\varepsilon,j})|\nabla\dot{\varphi}_{\varepsilon,j}|_{\omega_{\varepsilon,j}}^2\leq C|\nabla\dot{\varphi}_{\varepsilon,j}|_{\omega_{\varepsilon, j}}^2-|\nabla\overline{\nabla}\dot{\varphi}_{\varepsilon, j}|_{\omega_{\varepsilon, j}}^2-|\nabla\nabla\dot{\varphi}_{\varepsilon, j}|_{\omega_{\varepsilon, j}}^2+C.
\end{equation}
Set $H_{\varepsilon,j}=(t-\frac{\eta}{2})|\nabla\dot{\varphi}_{\varepsilon,j}|_{\omega_{\varepsilon,j}}^2+A\dot{\varphi}_{\varepsilon,j}^2$. Assume that the maximum point of $H_{\varepsilon,j}$ on $[\frac{\eta}{2}, T]$ is $(t_{0},x_{0})$ and $t_{0}>\frac{\eta}{2}$. Then we have
\begin{equation*}
\begin{split}
(\frac{d}{dt}-\Delta_{\varepsilon,j})H_{\varepsilon,t}
\leq &Ct|\nabla\dot{\varphi}_{\varepsilon,j}|^2+(1-2A)|\nabla\dot{\varphi}_{\varepsilon,j}|_{\omega_{\varepsilon,j}}^2+2A\dot{\varphi}_{\varepsilon,j}F'(\varphi_{\varepsilon,j}, z)\dot{\varphi}_{\varepsilon,j}+Ct\\
\leq &-|\nabla\dot{\varphi}_{\varepsilon,j}|_{\omega_{\varepsilon,j}}^2+C,\\
\end{split}
\end{equation*}
where $2A=CT+2$. By the maximum principle, we have
\begin{equation}
|\nabla\dot{\varphi}_{\varepsilon,j}|_{\omega_{\varepsilon,j}}^2\leq \frac{C}{t-\frac{\eta}{2}}.
\end{equation}
Straightforward calculations show that
\begin{equation*}
(\frac{d}{d t}-\Delta_{\varepsilon,j})\Delta_{\varepsilon,j} \dot{\varphi}_{\varepsilon,j}(t)\leq-|\nabla\overline{\nabla}\dot{\varphi}_{\varepsilon, j}|_{\omega_{\varepsilon, j}}^2+F'\Delta_{\varepsilon,j} \dot{\varphi}_{\varepsilon,j}(t)+C|\nabla\dot{\varphi}_{\varepsilon, j}|_{\omega_{\varepsilon, j}}+C.
\end{equation*}
Let $G_{\varepsilon,j}=(t-\frac{\eta}{2})^2(-\Delta_{\varepsilon, j}\dot{\varphi}_{\varepsilon, j})+2(t-\frac{\eta}{2})^2|\nabla\dot{\varphi}_{\varepsilon, j}|_{\omega_{\varepsilon, j}}^2$. Assume that the minmum point of $G_{\varepsilon,j}$ on $[\frac{\eta}{2}, T]$ is $(t_{0},x_{0})$, $t_{0}>\frac{\eta}{2}$ and $-\Delta_{\varepsilon,j}\dot{\varphi}(t_0, x_0)< 0$. Then we have
\begin{equation*}
\begin{split}
(\frac{d}{dt}-\Delta)G_{\varepsilon, j}\geq & (t-\frac{\eta}{2})(2+(t-\frac{\eta}{2})F')(-\Delta\dot{\varphi}_{\varepsilon, j})+|\nabla\overline{\nabla}\dot{\varphi}_{\varepsilon, j}|_{\omega_{\varepsilon, j}}^2(t-\frac{\eta}{2})^2-C\\
\geq &(t-\frac{\eta}{2})C(-\Delta_{\varepsilon, j}\dot{\varphi}_{\varepsilon, j})+\frac{(t-\frac{\eta}{2})^2}{n}(-\Delta_{\varepsilon, j}\dot{\varphi}_{\varepsilon, j})^2-C\\
=&\big(\frac{(t-\frac{\eta}{2})(-\Delta_{\varepsilon, j}\dot{\varphi}_{\varepsilon, j})}{\sqrt{n}}+\frac{\sqrt{n}C}{2}\big)^2-\frac{C^2n}{4}-C.
\end{split}
\end{equation*}
By the maximum principle, we get
\begin{equation}
(t-\frac{\eta}{2})^2(-\Delta_{\varepsilon, j}\dot{\varphi}_{\varepsilon, j})\geq -C.
\end{equation}
On the other hand, we set $\tilde{G}_{\varepsilon,j}=(t-\frac{\eta}{2})^2(-\Delta_{\varepsilon, j}\dot{\varphi}_{\varepsilon, j})+2(t-\frac{\eta}{2})^2|\nabla\dot{\varphi}_{\varepsilon, j}|_{\omega_{\varepsilon, j}}^2$. Assume that the maximum point of $\tilde{G}_{\varepsilon,j}$ on $[\frac{\eta}{2}, T]$ is $(t_{0},x_{0})$, $t_{0}>0$ and $-\Delta_{\varepsilon,j}\dot{\varphi}(t_0, x_0)>0$.
\begin{equation*}
\begin{split}
(\frac{d}{dt}-\Delta_{\varepsilon,j})\tilde{G}_{\varepsilon, j}\leq&(t-\frac{\eta}{2})(2+(t-\frac{\eta}{2})F')(-\Delta_{\varepsilon,j}\dot{\varphi}_{\varepsilon,j})
-2(t-\frac{\eta}{2})^2|\nabla\overline{\nabla}\dot{\varphi}_{\varepsilon,j}|_{\omega_{\varepsilon,j}}^2+C(t-\frac{\eta}{2})\\
\leq& C(t-\frac{\eta}{2})(-\Delta_{\varepsilon,j}\dot{\varphi}_{\varepsilon,j})-2(t-\frac{\eta}{2})^2\frac{(-\Delta_{\varepsilon,j}\dot{\varphi}_{\varepsilon,j})^2}{n}+C(t-\frac{\eta}{2}).
\end{split}
\end{equation*}
By the maximum principle, we get
\begin{equation}
(t-\frac{\eta}{2})^2(-\Delta_{\varepsilon,j}\dot{\varphi}_{\varepsilon,j})\leq C.
\end{equation}

At last, we prove the regularity along the conical parabolic complex Monge-Amp\`ere equation $(\ref{1})$.

\begin{thm} \label{2016011202}The solution $\varphi(t,\cdot)$ of the equation $(\ref{1})$ is $C^{2,\alpha,\beta}$ for any $\alpha\in(0,\min\{1,\frac{1}{\beta}-1\})$ when $t\in(0,T]$.
\end{thm}

{\bf Proof :}\ \ Fix $t>0$. We assume that $t\in[\eta,T]$ for some $0<\eta<T$. By $(\ref{222})$, we know that $|\Delta_{\omega_\beta}\varphi(t)|$ is uniformly bounded by some constant $C$ depending only on $\|\varphi_0\|_{L^\infty(M)}$, $n$, $\beta$, $\omega_{0}$, $\eta$, $F$ and $T$. By Proposition $\ref{21000000}$, we obtain that $|\nabla\varphi(t)|_{\omega(t)}$ is bounded. By the Lemma $4.6$ in \cite{JWLXZ}, Lemma \ref{2016080701}, Lemma $\ref{1.8.6}$ and $(\ref{222})$, we have
\begin{eqnarray*}\|\varphi\|_{C^1(M)}\leq C,\ \ \|\dot{\varphi}\|_{C^1(M)}\leq C,\ \ |\Delta_{\omega(t)}\dot{\varphi}|\leq C,\ \ \frac{1}{C}\omega_\beta\leq\omega_0+\sqrt{-1}\partial\bar{\partial}\varphi(t)\leq C\omega_\beta,
\end{eqnarray*}
where the $C^1$-norm and Laplacian are taken with respect to $\omega(t)$, and the constant $C$ depending only on $\|\varphi_0\|_{L^\infty(M)}$, $n$, $\beta$, $\omega_{0}$, $\eta$, $F$ and $T$. Since $\omega_\beta\geq\gamma\omega_0$ for some positive constant $\gamma$, we have
\begin{eqnarray}|\nabla F(\varphi(t,z),z)|_{\omega(t)}\leq C\ \ \ and\ \ \ n-C\gamma^{-1}n\leq\Delta_{\omega(t)}\varphi\leq n.
\end{eqnarray}
At the same time, there exists a constant $C$ such that
\begin{eqnarray}|\Delta_{\omega(t)}F(\varphi(t,z),z)|\leq C.
\end{eqnarray}
Now we consider the parabolic Monge-Amp\`ere equation
\begin{eqnarray}(\omega_0+\sqrt{-1}\partial\bar{\partial}\varphi)^n=e^{\dot{\varphi}-F(\varphi(t,z),z)}\frac{\omega_0^n}{|s|_h^{2(1-\beta)}}.
\end{eqnarray}
If we use Tian's theorem in \cite{Tian5} directly, we need prove that $\Delta_{\omega_\beta}\dot{\varphi}$ is bounded from below. But here we can only obtain the lower bound of $\Delta_{\omega(t)}\dot{\varphi}$, so we need a modification of Tian's theorem. Now replacing the base metric $\omega_{\beta}$ and the covariant derivatives with respect to $\omega_\beta$ in Tian's theorem \cite{Tian5} by $\omega(t)$ and ordinary partial derivatives respectively, we can conclude that for any $\alpha\in(0,\min\{1,\frac{1}{\beta}-1\})$,  $\sqrt{-1}\partial\bar{\partial}\varphi$ is $C^{\alpha}$-bounded with respect to $\omega(t)$ by following Tian's discussion. Since this process is very similar to the one in \cite{Tian5}, we omit the details here. Combining with the equivalence between $\omega(t)$ and $\omega_\beta$, we know that the unique solution to the conical parabolic complex Monge-Amp\`ere equation $(\ref{1})$ is $C^{2, \alpha,\beta}$ for any $\alpha\in(0,\min\{1,\frac{1}{\beta}-1\})$ when $t\in(0,T]$. \QEDB

\section{Proof of Theorem \ref{thm06}}
\setcounter{equation}{0}

In this section, combining the arguments in section $2$, $3$ with \cite{GSVT}, we prove Theorem \ref{thm06}. By Kolodziej's results \cite{K11,K}, we know that $\varphi_{0}$ is H\"older continuous on $M$ and there exists a continuous solution $\psi_{\varepsilon}$ to equation
\begin{eqnarray}
(\omega_{0}+\sqrt{-1}\partial\bar{\partial}\psi_{\varepsilon})^{n}=c_{\varepsilon}e^{-F(\varphi_{0},z)}\frac{\omega_{0}^{n}}{(\varepsilon^{2}+|s|_{h}^{2})^{1-\beta}},
\end{eqnarray}
where $\sup\limits_{M}(\varphi_0-\psi_{\varepsilon})=\sup\limits_{M}(\psi_{\varepsilon}-\varphi_{0})$, $c_{\varepsilon}$ is the normalized constant and it is uniformly bounded. Let $\{\varphi_{j}\}$ be a sequence of smooth functions such that
\begin{eqnarray}\lim\limits_{j\rightarrow\infty}\|\varphi_{j}-\varphi_{0}\|_{L^{\infty}(M)}=0.
\end{eqnarray}
By Yau's results \cite{STY}, there exists a smooth function $\psi_{\varepsilon,j}$ such that
\begin{eqnarray}
(\omega_{0}+\sqrt{-1}\partial\bar{\partial}\psi_{\varepsilon,j})^{n}=c_{\varepsilon,j}e^{-F(\varphi_{j}(z),z)}\frac{\omega_{0}^{n}}{(\varepsilon^{2}+|s|_{h}^{2})^{1-\beta}},
\end{eqnarray}
where $\sup\limits_{M}(\psi_{\varepsilon}-\psi_{\varepsilon,j})=\sup\limits_{M}(\psi_{\varepsilon,j}-\psi_{\varepsilon})$, $c_{\varepsilon,j}$ is the normalized constant and it is uniformly bounded.  By the stability theorem in \cite{K00} (see also \cite{SDZZ} or \cite{VGAZ1}),
\begin{eqnarray}\label{00000}
\lim\limits_{j\rightarrow\infty}\|\psi_{\varepsilon,j}-\psi_{\varepsilon}\|_{L^{\infty}(M)}&=&0,\\
\label{00000010}\lim\limits_{j\rightarrow\infty}\|\psi_{\varepsilon}-\varphi_{0}\|_{L^{\infty}(M)}&=&0.
\end{eqnarray}
\begin{pro} \label{2016081701} There exist uniform $T$ such that  the equation
\begin{equation}\label{CMAE5}
\left \{\begin{split}
&\frac{\partial \varphi_{\varepsilon,j}(t)}{\partial t}=\log\frac{(\omega_0+\sqrt{-1}\partial\bar{\partial}\varphi_{\varepsilon,j}(t))^{n}}{\omega_{0}^{n}}+F(\varphi_{\varepsilon, j}(t), z),\\
&\ \ \ \ \ \ \ \ \ \ \ \ \ \ +(1-\beta)\log(\varepsilon^2+|s|_h^2)-\log c_{\varepsilon,j}+\log c_{\varepsilon}\\
&\varphi_{\varepsilon,j}(0)=\psi_{\varepsilon,j}\\
\end{split}
\right.
\end{equation}
admits a unique smooth solution $\varphi_{\varepsilon,j}(t)$ on $[0,T]\times M$ for any $\varepsilon$ and $j$. The solutions $\{\varphi_{\varepsilon,j}(t)\}$ converge to the unique solution  $\varphi_\varepsilon(t)$ of equation
\begin{equation}\label{CMAE6}
\left \{\begin{split}
&\frac{\partial \varphi_{\varepsilon}(t)}{\partial t}=\log\frac{(\omega_0+\sqrt{-1}\partial\bar{\partial}\varphi_{\varepsilon}(t))^{n}}{\omega_{0}^{n}}+F(\varphi_{\varepsilon}(t), z)+(1-\beta)\log(\varepsilon^2+|s|_h^2),\\
&\varphi_{\varepsilon}(0)=\psi_{\varepsilon}\\
\end{split}
\right.
\end{equation}
in $L^{\infty}([0,T]\times M)$ and $C^{\infty}([\eta, T]\times M)$ sense for any $0<\eta<T$. Furthermore, $\varphi_{\varepsilon}(t)$ converge to $\varphi(t)$ which is the solution of equation $(\ref{1})$ with initial data $\varphi_{0}$ in $C^{\infty}_{loc}((0,T]\times (M\setminus D))$ sense.
\end{pro}
{\bf Proof:}\ \  The proof is similar as that in the section $2$ and $3$. We only prove the last result. By the arguments in section $3$, we conclude that $\varphi_{\varepsilon}(t)$ converges to a function $\tilde{\varphi}(t)$ in $C^{\infty}_{loc}((0,T]\times (M\setminus D))$. We prove $\tilde{\varphi}(t)=\varphi(t)$ by using the uniqueness Theorem $\ref{228}$.  By checking the steps in section $3$, we know that all the arguments are valid except Proposition $\ref{218}$. So we only need prove inequality $(\ref{201808})$. By the similar arguments as that in Proposition \ref{21000000}, for any $\varepsilon_1<\varepsilon_2$, $(t,x)\in [0,T]\times M$ and $\delta>0$, we have
\begin{eqnarray}\varphi_{\varepsilon_{1}}(t,x)-\varphi_{\varepsilon_{2}}(t,x)\leq e^{KT}\|\psi_{\varepsilon_{1}}-\psi_{\varepsilon_{2}}\|_{L^{\infty}(M)}+Te^{K T}\delta.
\end{eqnarray}
For any $t>0$, we let $\delta\rightarrow0$ and then $\varepsilon_{1}\rightarrow0$,
\begin{eqnarray*}\tilde{\varphi}(t,x)-\varphi_{0}&\leq& \varphi_{\varepsilon_{2}}(t,x)+e^{KT}\|\varphi_{0}-\psi_{\varepsilon_{2}}\|_{L^{\infty}(M)}-\varphi_{0}\\
&\leq&|\varphi_{\varepsilon_2}(t,z)-\varphi_{\varepsilon_2,j}(t,z)|
+|\varphi_{\varepsilon_2,j}(t,z)-\psi_{\varepsilon_{2},j}(z)|\\
&\ &+|\psi_{\varepsilon_{2},j}(z)-\psi_{\varepsilon_{2}}|+|\psi_{\varepsilon_{2}}-\varphi_{0}(z)|+e^{KT}\|\varphi_{0}-\psi_{\varepsilon_{2}}\|_{L^{\infty}(M)}.
\end{eqnarray*}
By choosing suitable $\varepsilon_{2}$, $j$ and $\delta_{1}$,  we have
\begin{eqnarray}\tilde{\varphi}(t,z)-\varphi_{0}(z)<\epsilon\ \ for\ any\ (t,x)\in(0,\delta_1]\times M.
\end{eqnarray}
Then by constructing auxiliary function as that in Proposition $\ref{101}$,
\begin{eqnarray}\lim\limits_{t\rightarrow0}\|\tilde{\varphi}(t)-\varphi_{0}\|_{L^{\infty}(M)}=0.
\end{eqnarray}
By the uniqueness Theorem $\ref{228}$, we have $\tilde{\varphi}(t)=\varphi(t)$.\QEDB

{\bf Proof of Theorem $\ref{thm06}$:}\ \  Differentiating the equation $(\ref{CMAE5})$ we obtain
\begin{eqnarray}(\frac{d}{dt}-\Delta_{\varepsilon,j})\dot{\varphi}_{\varepsilon,j}(t)=F'(\varphi_{\varepsilon,j}(t),z)\dot{\varphi}_{\varepsilon,j}(t).
\end{eqnarray}
Since $F'(\varphi_{\varepsilon,j}(t), z)$ is bounded as long as $\varphi_{\varepsilon,j}(t)$ is
bounded, from the maximum principle we get
\begin{eqnarray}\sup\limits_{M}|\dot{\varphi}_{\varepsilon,j}(t)|\leq e^{Kt}\sup\limits_{M}|\dot{\varphi}_{\varepsilon,j}(0)| \ for\ any\ t\in [0,T],
\end{eqnarray}
where $K$ depends on $\beta$, $n$, $\omega_{0}$, $F$, $T$ and $\|\varphi\|_{L^{\infty}}$. Since
\begin{eqnarray}\dot{\varphi}_{\varepsilon,j}(0)=F(\psi_{\varepsilon, j}, z)-F(\varphi_{j}(z), z)+\log c_{\varepsilon},
\end{eqnarray}
for any $(t,z)\in(0,T]\times(M\setminus D)$, we have
\begin{eqnarray}|\dot{\varphi}_{\varepsilon,j}(t,z)|\leq e^{KT}(|F(\psi_{\varepsilon, j}, z)-F(\varphi_{j}(z), z)|+|\log c_{\varepsilon}|).
\end{eqnarray}
Let $j\rightarrow \infty$ and then $\varepsilon\rightarrow 0$, by Proposition $\ref{2016081701}$, we obtain
\begin{eqnarray}\dot{\varphi}(t,z)=0\ on\ (0,T]\times (M\setminus D).
\end{eqnarray}
Since $\varphi(t,z)$ is continuous on $[0, T]\times M$, it follows that $\varphi(t,z)=\varphi_{0}$ for all $t$. By Theorem $\ref{thm04}$, $\varphi_{0}$ is $C^{2,\alpha,\beta}$ for any $\alpha\in(0,\min\{1,\frac{1}{\beta}-1\})$.\QEDB

\hspace{1.4cm}

\end{document}